\documentclass[11pt,a4paper]{article}
\usepackage[T1]{fontenc}
\usepackage[latin1]{inputenc}
\usepackage{amsfonts}
\usepackage{latexsym}
\usepackage{graphicx}
\usepackage{amssymb}
\usepackage{amsfonts}
\usepackage{amsmath}
\usepackage{amsthm}
\usepackage[english]{babel}
\usepackage{epsfig}
\usepackage{verbatim}
\newcommand{\pr}{\mathbf P}
\newcommand{\ew}{{\mathbf E}}
\newcommand{\cov}{{\mathbf{Cov}}}
\newcommand{\var}{{\mathbf{Var}}}
\newcommand{\vr}{{\varrho}}
\newcommand{\R}{{\mathbb R}} 

\renewcommand{\a}{\alpha}
\newcommand{\g}{\gamma}
\newcommand{\e}{\varepsilon}
\renewcommand{\t}{\theta}
\renewcommand{\l}{\lambda}

\newcommand{\pa}{\partial}

\renewcommand{\d}{\delta}
\newcommand{\D}{\Delta}

\newcommand{\Su}{\mathop{{\sum}^{\ne}}}
\newcommand{\Plongn}{\stackrel{{\mathbf P}}{\longn}}

\renewcommand{\O}{\mathbf o}
\newcommand{\1}{\mathbf 1}

\newcommand{\mtl}{\mathop{\times}\limits}
\newcommand{\longn}{\mathop{\longrightarrow}\limits_{n \to \infty}}

\newcommand{\dlongn}{\stackrel{d}{\longn}}
\newcommand{\longP}{\stackrel{{\pr}-a.s.}{\longn}}
\newcommand{\be}{\begin{equation}}
\newcommand{\ee}{\end{equation}}
\newcommand{\bea}{\begin{eqnarray}}
\newcommand{\eea}{\end{eqnarray}}
\newcommand{\n}{\nonumber}
\newcommand{\sn}{\smallskip\noindent}
\newcommand{\bn}{\bigskip\noindent}
\newcommand{\mn}{\medskip\noindent}

\begin{document}




\title{\vspace*{-1 cm} 
\bf Asymptotic Goodness-of-Fit Tests for Point Processes \\  Based on Scaled Empirical K-Functions}

\author{{\bf Lothar Heinrich}$^1$}

\date{}
\maketitle

\begin{abstract}
\noindent
{\small We study sequences of scaled edge-corrected empirical (generalized) $K$-functions 
(modifying Ripley's $K$-function) each of them constructed from a single observation of a 
$d$-dimensional fourth-order stationary point process in a sampling window $W_n$ which 
grows together with some scaling rate unboundedly as $n \to \infty$. Under some natural 
assumptions it is shown that the normalized difference between scaled empirical and scaled 
theoretical $K$-function converges weakly to a mean zero Gaussian process with simple 
covariance function. This result suggests discrepancy measures between empirical and theoretical  
$K$-function with known limit distribution which allow to perform goodness-of-fit tests 
for checking a hypothesized point process based only on its intensity and (generalized) $K$-function. 
Similar test statistics are derived for testing the hypothesis that two independent point processes 
in $W_n$ have the same distribution without explicit knowledge of their intensities and $K$-functions.} 


\vspace{0.3cm}
 
\vspace{0.3cm}
\end{abstract}


{\vspace*{-1 cm}
\stepcounter{footnote}\footnotetext{University of Augsburg, Institute of Mathematics, Universit\"atsstr. 14,
86135 Augsburg, Germany;\\ \phantom{000} E-mail: heinrich@math.uni-augsburg.de}


\vspace{8mm}
\section{Introduction} 

\bn
In many fields of application statisticians are faced with irregular point
patterns  or  point-like objects which are randomly distributed  in large 
planar or spatial domains. Random point processes provide appropriate 
models to describe such phenomena. It is often assumed and in practice at 
least approximately justifiable that the distribution of the point pattern 
under consideration or at least its  moment measures up to certain order
are invariant under translations. In this paper all point processes 
and random variables are defined over a common probability space $[\Omega, {\mathcal A}, \pr]$
and $\ew\,$, $\var\,$, $\cov$ denote expectation, variance, covariance w.r.t. $\pr\,$.

\medskip
The aim of the paper is to establish asymptotic goodness-of-fit tests for checking point 
process hypotheses provided the hypothesized $d$-dimensional simple point process (short: PP)
$N = \sum_{i\ge 1}\delta_{X_i}$ with distribution $P$ (short: $N \sim P$) 
is fourth-order stationary  with known intensity $\l := \ew N([0,1)^d) > 0$ and 
a known generalized $K$-function $0 \le r \mapsto K_B(r) := \l^{-1}\,\a_{red}^{(2)}(r B)\,$,  
where $\a_{red}^{(2)}(\cdot)$ denotes the reduced second factorial moment measure which always
exists when $N$ is at least second-order (or weakly) stationary, see Daley and Vere-Jones (1988) 
for mathematical background of PP theory. Here, $B$ is assumed to be an $\O$-symmetric, convex, 
compact set containing the origin $\O \in \R^d$ as inner point. 
Equivalently, $B$ can be considered as closed unit ball w.r.t. a unique norm $\|\cdot\|_B$ 
defined by $\|x\|_B := \inf\{r > 0: x \in r B\}$ for $x \in \R^d$. In case of the Euclidean
norm  $\|\cdot\|$ with unit ball $B_e$ the  subscript $B$ is omitted and $K(r)$ coincides 
with Ripley's  $K$-function as introduced in Ripley (1976) for motion-invariant PPes. 
Our approach is based on a single observation of $N$ in a large sampling window $W_n$, where
throughout this paper $(W_n)$  forms a {\em convex averaging sequence} (short: CAS), i.e. 
$W_n \subseteq W_{n+1}$  and $W_n $ is convex, compact in $\R^d$ for all $n \ge 1$ such that 
$\vr(W_n):=\sup\{r>0 : x+r\,B_e \subseteq W_n$ for some $x \in W_n\} \longn \infty$.

\medskip
Instead to impose one of the usual strong mixing conditions on $N \sim P$ we will assume that the    
the second, third and fourth reduced cumulant measure have bounded total variation and that the ratio
$(N(W_n)-\l\,|W_n|)/\sqrt{|W_n|}$ is asymptotically normally distributed (as $n \to \infty$)
with an {\em asymptotic variance}  $\sigma^2 := \lim\limits_{n\to \infty}\var(N(W_n))/|W_n|\,$. 
Here and in what follows, $|\cdot|$ stands for the Lebesgue measure or volume in $\R^d$ as well as 
for the absolute value of a real number.  
 
A frequently used unbiased (edge-corrected) estimator of $\l^2\,K_B(r)$ (so far for $B=B_e$) is the so-called {\em Ohser-Stoyan-} 
or {\em Horvitz-Thompson-type estimator} defined by

\be
(\widehat{\l^2 K_B})_n(r) := \Su_{i,j \ge 1} \frac{\1_{W_n}(X_i)\,\1_{W_n}(X_j)\;\1_{r\,B}( X_j - X_i )}{|(W_n - X_i) \cap (W_n - X_j)|}\,,
\label{KBest}\ee

\noindent
see Ohser and Stoyan (1981), Stein (1993), Diggle (2003), Chiu et al. (2013), Chapt. 4.7.4, or Illian et al. (2008), Chapt. 4.3.3, 
for further details and modifications. On the r.h.s of (\ref{KBest}) the sum $\Su$ runs over pairs of distinct indices and 
the indicator function ${\bf 1}_{r\,B}( X_j - X_i )$  can be replaced by 
$\1_{[0,r]}(\| X_j - X_i \|_B)$. It is plausible that the use of several  norms $\|\cdot\|_B$ provide more information 
on the PP to be examined. In this connection $B$ plays a similar role like {\em structuring elements} in image analysis, see Chiu et al.(2013). 
It should be mentioned that there exist a number of further estimators based on second-order characteristics to obtain information on stationary
point patterns (orientation, isotropy etc.), see e.g. Guan et al. (2006), Guan and Sherman (2007), Illian et al. (2008), Chapt. 4.5.3.
Baddeley et al. (2000) introduced an inhomogeneous $K$-function and its empirical counterpart for intensity-reweighted stationary PPes, see also
Gaetan and Guyon (2010), Chapt. 3.5.3, Adelfio and Schoenberg (2009) and Zhao and Wang (2010) for asymptotic properties.  

\medskip
In the most important case when $N \sim P$ is a stationary Poisson process the 
weak convergence of the empirical process $\sqrt{|W_n|}\,\big(\,(\widehat{\l^2 K})_n(r) - \l^2 |B_e| r^d\,\big)$ 
in the Skorohod space $D[0,R]$, see Billingsley (1968), to a Gaussian process has been proved in Heinrich (1991).
The corresponding functional central limit (short: CLT) for the empirical process (\ref{KBest}) can be show in the same way, 
see Heinrich (2013). The known limit distribution of $\sqrt{|W_n|}\,\max_{0\le r \le R}\big|\,(\widehat{\l^2 K_B})_n(r) - 
\l^2 |B| r^d\,\big|$ (in particular when $\l^2$ is replaced by an unbiased estimator) allows to perform  Kolmogorov-Smirnov-type
as well as Cram\'{e}r-von Mises-type tests to check {\em complete spatial randomness (CSR)} of a given point pattern in $W_n$.  Multiparameter versions of this and related 
asymptotic tests are derived in Heinrich (2015). Note that testing for CSR of a transformed PP can provide information on the original PP, 
see Schoenberg (1999). There exist further (non-asymptotic) tests for CSR based on Ripley's K-function using 
different estimators supported by data analytic methods and simulation techniques, see e.g. Ho and Chiu (2010), Marcon et al. (2013), 
Wiegand et al. (2016) and further references therein.

\medskip  
Asymptotic normality of the family of random variables $\sqrt{|W_n|}\,\big((\widehat{\l^2 K_B})_n(r) - \l^2 K_B(r)\big) $ 
can be shown due to a more general CLT  proved first by Jolivet (1981), see also Karr (1987) and Ki\^{e}u and Mora (1999), if the underlying 
PP $N \sim P$ is  Brillinger mixing, see Section 2. The main obstacle to construct a goodness-of-fit test as in the Poisson case 
is the rather complicated covariance function of the Gaussian limit process depending on $\l$ and integrals w.r.t. $\g_{red}^{(k)}$, $k=2,3,4\,$, 
see (\ref{gred}). A  closed-term formula of this covariance function is given at the end of Appendix B, whereas a more explicit form of this function is known e.g. 
for $\a$-determinental PPes, see Heinrich (2016). An  expression for the asymptotic variance of $\sqrt{|W_n|}\,(\widehat{\l^2 K_B})_n(r)$ for $r > 0$ is given in (\ref{asvar}).

\medskip
To realize the announced  program of testing hypotheses on the distribution of non-Poisson PPes via their generalized $K$-functions $K_B(r)$ 
we allow that  the range of $r$  grows  unboundedly with a scaling rate $c_n^{\a}\,$, where $c_n := |W_n|^{1/d}$ and $0 < \a < 1$.  
More precisely, our main issue is to study the asymptotic behaviour of the sequence of the scaled (centered and normalized) empirical processes 

\be 
\D^{(\a)}_{B,n}(r) := |W_n|^{1/2 - \a}\,\big(\,(\widehat{\l^2 K_B})_n(c_n^{\a}\,r) - \l^2\,K_B(c_n^{\a}\,r)\,\big)\quad\mbox{for}\quad 0 \le r \le R
\;\;,\;\;0 \le \a < 1
\label{Dna}\ee 

\mn
as $n \to \infty\,$, where $R$ is a freely chosen positive real number. 
Under mild (mixing) assumptions we shall show in Section 3 that the asymptotic variance and covariance of (\ref{Dna}) does not 
depend on $\g_{red}^{(3)}$ and $\g_{red}^{(4)}$.  Furthermore, (\ref{Dna}) possesses a Gaussian limit for fixed $r > 0$ and in 
$D[0,R]$) for all $0 < \a < 1$, but only for $\a \le 1/2$ this results are relevant for our tests. The latter restriction is 
meaningful because $\|\g_{red}^{(2)}\|^{}_{var} < \infty$ implies for $\a > 1/2$ that

\be 
\sup\limits_{0 \le r \le R} |W_n|^{1/2-\a}\,\big|\,\l^2 K_B(c_n^{\a}\,r) - \l^2\,| c_n^{\a}\,r\,B |\,\big|  = 
\l\,|W_n|^{1/2-\a}\,\sup\limits_{0 \le r \le R}\big|\,\g_{red}^{(2)}(c_n^{\a}\,r\,B)\,\big| \longn 0\,.
\label{remainder}\ee

\mn
This means that we can always take $K_B(r) = |r\,B| = |B|\,r^d$ in (\ref{Dna}) which in turn prevents to distinguish between Poisson and
non-Poisson PPes if $\a > 1/2\,$.  On the other hand, such a separation is possible for $\a \in (0,1/2]\,$, where the limit distribution of 
(\ref{Dna}) depends on the intensity $\l$ and asymptotic variance $\sigma^2$ of the PP $N \sim P$.  
In particular for $\a = 1/2$, under comparatively mild assumptions on $N \sim P$ we can show that  

\be 
\frac{1}{\widehat{\l_n}\,\sqrt{(\widehat{\sigma^2})_n}} \;\max\limits_{0 \le r \le \sqrt{c_n}\,R}
\big|\,(\widehat{\l^2 K_B})_n(r) - \l^2\,K_B(r)\,\big| \dlongn 2\,R^d\,|B|\;|{\mathcal N}(0,1)|\;,
\label{KBclt0.5}\ee

\mn
where $(\widehat{\sigma^2})_n$ and $\widehat{\l_n}$  are defined in (\ref{sigest}) resp. (\ref{l2est}),  ${\mathcal N}(0,\tau^2)$ denotes 
a mean zero Gaussian random variable with variance $\tau^2$ and $\dlongn$ indicates {\em convergence in distribution} or {\em weak convergence}.

\medskip
More general limit theorems are presented, proved and applied  in the Sections 3 - 5 to construct the above-announced goodness-of-fit tests.
This includes tests of the null hypothesis that the distributions $P_a$ and $P_b$ of two independent stationary PPes $N_a$ and $N_b$ 
coincide based on checking the coincidence of their $K$-functions $K_{a,B}$ and $K_{b,B}$ and their intensities $\l_a$ and $\l_b\,$. 

\vspace{5mm}
\section{Definitions and Preliminary  Results}

\mn 
A simple point process $N = \sum_{i\ge 1} \d_{X_i}$ on $\R^d$ 
is said to be $k${\em th-order stationary} for some $k \ge 1$ if 
${\ew}N^k([0,1)^d) < \infty$ and, for $j=1,...,k\,$,
\be 
\a^{(j)}(\mtl_{i=1}^j (A_j+x)) = \a^{(j)}(\mtl_{i=1}^j A_i)\quad\mbox{for all}\quad
x\in R^d\quad\mbox{and any bounded}\quad A_1,...,A_j \in {\mathcal B}^d\,.
\n\ee
Here $\a^{(k)}(\cdot)$ denotes the  $k$-{\em th factorial moment measure}
which is defined on $[\R^{dk},{\mathcal B}^{dk}]$  by
\be
\a^{(k)}(\mtl_{i=1}^k A_i) = \ew \Su_{i_1,...,i_k \ge 1} \1_{A_1}(X_{i_1})\cdots \1_{A_k}(X_{i_k})\;, 
\n\ee
where the sum $\Su$ stretches over $k$-tuples of  pairwise distinct indices.

\bigskip
The invariance of $\a^{(k)}(\cdot)$ against diagonal-shifts (if $k\ge 2$) allows the disintegration 
w.r.t. the Lebesgue measure $|\cdot |$ (multiplied additionally by $\l$) and defines 
the {\em reduced} $k$-{\em th factorial moment measure} $\a_{red}^{(k)}(\cdot)$ 
(on $[\R^{d(k-1)},{\mathcal B}^{d(k-1)}]$) satisfying
\be 
\a^{(k)}(\mtl_{i=1}^k A_i) = \l \int_{A_k} \a_{red}^{(k)}(\mtl_{i=1}^{k-1}(A_i-x))\,{\rm d}x 
\quad\mbox{for any bounded}\quad A_1,...,A_k \in {\mathcal B}^d\,,
\n\ee

\mn
where, for bounded $W, A_2,...,A_k \in {\mathcal B}^d$ with $|W| > 0$, the measure $\a_{red}^{(k)}(\cdot)$ can be expressed by
\be 
\a_{red}^{(k)}(\mtl_{i=2}^k A_i) = \frac{1}{\l\,|W|}\,\ew \Su_{i_1,...,i_k \ge 1} \1_{W}(X_{i_1})\,\1_{A_2}(X_{i_2}- X_{i_1})\cdots \1_{A_k}(X_{i_k}- X_{i_1})\,.
\n\ee

Next we recall formal definition of $k$-{\em th factorial cumulant measure} $\g^{(k)}(\cdot)$ 
which is a signed measure on $[\R^{dk},{\mathcal B}^{dk}]$ defined by 
\be 
\g^{(k)}(\mtl_{i=1}^k A_i) = \sum_{j=1}^k (-1)^{j-1}\,(j-1)!\,\sum_{K_1\cup\cdots\cup K_j
=\{1,..,k\}}\prod_{i=1}^j \a^{(\#K_i)}(\mtl\limits_{k_i\in K_i} A_{k_i})\;,
\label{ga}
\ee

\mn
where the inner sum is taken over all partitions of the set $\{1,...,k\}$
into $j$ disjoint non-empty subsets $K_1,...,K_j\,$, see Daley and Vere-Jones (1988), Heinrich (2013).

A corresponding representation formula of $\a^{(k)}(\cdot)$ in terms of factorial cumulant measures 
can be derived  from (\ref{ga}) by induction,
\be 
\a^{(k)}(\mtl_{i=1}^k A_i) = \sum_{j=1}^k\sum_{K_1\cup\cdots\cup K_j
=\{1,..,k\}}\prod_{i=1}^j \g^{(\#K_i)}(\mtl\limits_{k_i\in K_i} A_{k_i})\,.
\label{ag}\ee

\medskip
Since $k$th-order stationarity implies diagonal-shift invariance of the locally finite signed 
measures $\g^{(j)}(\cdot)$ for $j=1,..,k\,$, the above
disintegration also applies for factorial cumulant measures. In this way a 
$k$th-order stationary PP is accompanied by reduced factorial cumulant 
measures $\g_{red}^{(j)}(\cdot)$ of order $j=2,...,k$, defined by 

\be 
\g^{(j)}(\mtl_{i=1}^j A_i) = \l\,\int_{A_j} \g_{red}^{(j)}(\mtl_{i=1}^{j-1}(A_i-x)){\mathrm d}x\,.
\label{gred}\ee

For $j=1$ it makes sense to put  $\g_{red}^{(1)}(\cdot) = \a_{red}^{(1)}(\cdot) \equiv 1$, whereas 
the signed measure $\g_{red}^{(2)}(\cdot) = \a_{red}^{(2)}(\cdot) - \l\,|\cdot|$ ( called 
{\em reduced covariance measure}) reflects the dependence between  two distinct atoms of 
 $N$ and $\g_{red}^{(j)}(\cdot)$ describes mutual dependences between $j$ distinct atoms.

If both the positive and negative part 
$(\g_{red}^{(k)})^{\pm}(\cdot)$ of  Jordan's decomposition 
of $\g_{red}^{(k)}(\cdot)$ are finite on $\R^{d(k-1)}\,$, then 
\be 
\|\g_{red}^{(k)} \|_{var} := (\g_{red}^{(k)})^+(\R^{d(k-1)}) + (\g_{red}^{(k)})^-(\R^{d(k-1)})
\n\ee

\mn
is called the {\em total variation} of $\g_{red}^{(k)}(\cdot)$. A $k$th-order stationary PP $N = \sum_{i\ge 1} 
\d_{X_i}$ on $\R^d$  is said to be $B_k$-{\em mixing}  if $\|\g_{red}^{(j)}\|_{var} < \infty$ for $j=2,\ldots,k$.  
$B_k$-mixing for all $k \ge 2$ is called {\em Brillinger-mixing}, see e.g. Jolivet (1981), Karr (1987), Heinrich (2016), 
Biscio and Lanvancier (2016),(2017) for numerous applications in PP statistics. It should be mentioned that 
$B_k$-mixing of $N \sim P$ even for all $k \ge 2$ does not imply the uniqueness of the stationary distribution $P$, 
see also Baddeley and Silverman (1984) for $k=2\,$.  

\bigskip 
Next, we state some technical results which will be needed to prove the below Theorems 1-5. 

\bn
\textbf{Lemma 1} For any sequence  $r_n$ of positive numbers 
satisfying $r_n/\vr(W_n) \longn 0\,$, we have

\be 
\frac{1}{|W_n|}\,\sup_{x \in r_n B}\,\bigl(\,|W_n | - |W_n \cap (W_n - x)|\,\bigr) = 
{\mathcal O}\Big(\frac{r_n}{\vr(W_n)}\Big)\quad{as}\quad n \to \infty\,.
\n\ee

\mn
{\em Proof of Lemma 1} \ Since $B = \{x \in \R^d : \|x\|_B \le 1\}$ is bounded there exist some $\kappa > 0$ 
(depending on $B$) such that $B \subseteq \kappa B_e\,$. By definition of Minkowski-subtraction $\ominus$, 
see Chiu et al., p. 6, it follows that  

\be 
W_n \ominus (r\,\kappa\,B_e) \subseteq  W_n \ominus (r B) = \bigcap_{y \in r B} (W_n - y) 
\subseteq   W_n \cap (W_n - x)\quad\mbox{for any}\quad x \in r B\,.
\n\ee

\mn
Applying the first inequality (\ref{inequ}) of Lemma 5 to $K = W_n\,$ we see that

\be 
\sup\limits_{x \in r B}\Bigl( 1 - \frac{|W_n \cap (W_n - x)|}{|W_n|} \Bigr) 
\le 1 - \frac{|W_n \ominus (r\,\kappa\,B_e)|}{|W_n|} 
\le \frac{d\,\kappa\,r}{\vr(W_n)}\quad\mbox{for any}\quad r \in [0,\vr(W_n)/\kappa]\,.
\n\ee

\mn
Setting $r = r_n$ and letting $n \rightarrow \infty$ complete the proof of Lemma 1. $\Box$ 
   
\bigskip
The following Lemma 2 turns out crucial to prove the Theorems 1 and 2 in Section 3.

\bn
\textbf{Lemma 2} \ Let  $N = \sum_{i\ge 1} \d_{X_i}$ be a $B_4$-mixing PP
on $\R^d\,$. Further, assume that $c_n^{\a}/\vr(W_n) \longn 0$ for some $ 0 < \a < 1\,$. 
Then

\be
\lim_{n \to \infty}\,\ew \Bigl(\,\D^{(\a)}_{B,n}(r)
\,-\,2\,\frac{\a_{red}^{(2)}(c^{\a}_n\,r\,B)}{|W_n|^{\a}}\,\frac{N(W_n) - \l\,|W_n|}{\sqrt{|W_n|}}\,\Bigr)^2 \, = \,0 \quad\mbox{for all}\quad r \ge0\,.
\n\ee 

\bigskip 
The next result becomes relevant in proving the tightness of the sequence $(\D^{(\a)}_{B,n})_n\,$.

\bn
\textbf{Lemma 3} \ Let $N = \sum_{i\ge 1} \d_{X_i}$ be a $B_4$-mixing PP on $\R^d\,$. Further, 
let there exist a uniformly bounded Lebesgue density of $\a_{red}^{(2)}(\cdot)$, i.e. there exist 
a constant $a_0 > 0\,$ such that
\be
\a_{red}^{(2)}(A) \le a_0\,| A |\quad\mbox{for all bounded }\quad A \in {\mathcal B}^d \;.
\label{density}\ee 

\noindent
Then, for  $ 0 \le \a < 1$, $0 \le s < t \le R$ and all $n \ge n_0(R)\,$,

\be
\ew \bigl(\,\D_{B,n}^{(\a)}(t) - \D_{B,n}^{(\a)}(s)\,\bigr)^2\, \le \,a_1\,|(t B) \setminus (s B)|^2
+ a_2 \, |(t B) \setminus (s B)|\,|W_n|^{-\a} + a_3\,|W_n|^{-2\a}\,,
\label{lemma2}\ee

\sn
where the constants $a_1, a_2, a_3$ depend on $a_0,\l,\|\g_{red}^{(2)}\|_{var},\|\g_{red}^{(3)}\|_{var},\|\g_{red}^{(4)}\|_{var}\,$.

\bn
\textbf{Remark 1} \ If, in addition to (\ref{density}), the Lebesgue densities $\kappa^{(3)}$ and $\kappa^{(4)}$ of $\g_{red}^{(3)}(\cdot)$ 
resp. $\g_{red}^{(4)}(\cdot)$ exist such that, for all $x,y \in \R^d\,$,

\be
| \kappa^{(3)}(x,y) | \le a_{01}\quad\mbox{and}\quad \int_{\R^d}\int_{\R^d}\Big(\,\int_{\R^d}| \kappa^{(4)}(x,y+z,z) |\,
{\rm d}z\,\Big)^2 {\rm d}y{\rm d}x \le a_{02}\;,
\label{density34}\ee

\mn
then, for  $0 \le \a < 1$, $0 \le s < t \le R$ and all $n \ge n_0(R)\,$,

\be 
\ew \bigl(\,\D_{B,n}^{(\a)}(t) - \D_{B,n}^{(\a)}(s)\,\bigr)^2\, \le \,a_1\,|(t B) \setminus (s B)|^2
+ a_2 \,|(t B) \setminus (s B)|\,|W_n|^{-\a}\;,
\label{rem1}\ee

\sn
where the constants $a_0, a_1, a_2$ may differ from those in (\ref{density}) and (\ref{lemma2}), respectively. 

\bn
The proofs of Lemma 2 and  Lemma 3 as well as of (\ref{rem1}) are postponed to Section 5. 

\bigskip
Finally, we introduce an estimator of the asymptotic variance $\sigma^2 = \lim_{n\to\infty}\var(N(W_n))/|W_n|$ 
$= \l\,( 1 + \g^{(2)}_{red}(\R^d))$ which generalizes a corresponding estimator for cubic windows $W_n = [-n,n]^d$ 
defined and studied in Heinrich (1994), Heinrich and Proke\v{s}ov\'{a} (2010).

\bigskip
Let $w : \R^d \mapsto [0,\infty)$ be a Borel-measurable, bounded, symmetric function with bounded support satisfying 
$\lim_{x \to \O} w(x) = w(\O) = 1\,$, e.g. $w(x) = \1_B(x)$ or $w(x) = \bigl(\,1 - \|x \|\,\bigr)\,\1_B(x)\,$.

\be 
(\widehat{\sigma^2})_n := \widehat{\l}_n  + \Su_{i,j \ge 1}\frac{\1_{W_n}(X_i)\; \1_{W_n}(X_j)
}{|(W_n - X_i) \cap (W_n - X_j)|} \,w\Bigl(\frac{X_j - X_i}{b_n\,c_n}\Bigr)
- (\widehat {\l^2})_n\;(b_n\,c_n)^d \int_{\R^d} w(x) {\rm d}x\,,
\label{sigest}\ee
where $(b_n)$ is a positive null sequence of bandwidths such that $c_n\,b_n \longn \infty\,$, and

\be
\widehat{\l}_n\,:=\,\frac{N(W_n)}{|W_n|}\quad 
\mbox{and}\quad (\widehat {\l^2})_n\,:=\,\frac{(N(W_n) - 1)\,N(W_n)}{|W_n|^2}\,.
\label{l2est}\ee

The  asymptotic properties of $(\widehat{\sigma^2})_n$ summarized in

\bn
\textbf{Lemma 4} \ Let $(W_n)$ be a CAS in $\R^d$ with inball radius $\vr(W_n)$ and 
$N = \sum_{i \ge 1} \d_{X_i}$ be a $B_2$-mixing PP on $\R^d$. Then $(\widehat{\sigma^2})_n$
is asymptotically unbiased for $\sigma^2$, i.e.
\be 
\ew (\widehat{\sigma^2})_n \longn \sigma^2\,.
\n\ee

\mn
If, additionally, the PP $N$ is $B_4$-mixing and 
 $c_n b_n^2 \longn 0\;$ such that $c_n b_n/\vr(W_n) \longn 0\,$, 
then  $(\widehat{\sigma^2})_n$ is mean square consistent, i.e.
\be 
\ew \bigl(\,(\widehat {\sigma^2})_n - \sigma^2\,\bigr)^2 \longn 0\;.
\n\ee

\bn
{\em Proof of Lemma 4} \ In view of the definitions of factorial moment and cumulant measures and their reduced versions
it is easily seen after some rearrangements that

\be 
\ew (\widehat {\sigma^2})_n = \l + \l\,\int_{\R^d} w\Bigl(\frac{y}{b_n\,c_n}\Bigr)\,\g_{red}^{(2)}({\mathrm d}y)
- \frac{b_n^d}{|W_n|}\,\g^{(2)}(W_n\times W_n)\,\int_{\R^d} w(x)\,{\mathrm d}x\,. 
\n\ee

Since $\|\g_{red}^{(2)} \|_{var} < \infty$, $b_n \longn 0$ and $b_n\,c_n \longn \infty$, the properties of the function $w(\cdot)$ 
imply the limit  ${\ew}(\widehat {\sigma^2})_n\longn \l\,(1 + \g_{red}^{(2)}(\R^d))$. This proves the first assertion of Lemma 4.

The proof of the second assertion is somewhat more complicated because integrals w.r.t. third- and fourth-order factorial moment and 
cumulant measures have to be estimated. However, apart from obvious changes, this proof equals almost verbatim that of Theorem 1 in 
Heinrich and Proke\v{s}ov\'{a} (2010) (which states Lemma 4 for $W_n=[-n,n]^d$). 
Therefore we omit the details and refer the reader to Heinrich and Proke\v{s}ov\'{a} (2010).  $\Box$

\vspace{5mm}

\section{Limit Theorems for Scaled Empirical K-Functions}
 
\bn
\textbf{Theorem 1} \ Let  $N = \sum_{i\ge 1} \d_{X_i}$ be a $B_4$-mixing PP on $\R^d$ satisfying 

\be
\sqrt{|W_n|}\,(\,\widehat{\l}_n  - \l\,) = \frac{N(W_n) - \l \,|W_n|}{\sqrt{|W_n|}} \dlongn {\mathcal N}(0,\sigma^2)
\quad\mbox{with}\quad \sigma^2 = \l\,\bigl( 1 + \g^{(2)}_{red}(\R^d)\bigr)\,.
\label{cltN}\ee

\bn
If $c_n^{\a}/\vr(W_n) \longn 0$  for some $\a \in (0,1)\,$, then

\be 
\D_{B,n}^{(\a)}(r) \dlongn G_B(r) := 2\,\l\,|B| \,r^d\,{\mathcal N}(0,\sigma^2)\quad\mbox{for any fixed}\quad r \ge 0
\label{dDna}\ee
and
\be
(\D_{B,n}^{(\a)}(r_i))_{i=1}^k \dlongn (G_B(r_i))_{i=1}^k\quad\mbox{for any fixed $k$-tuple}\;\;r_1,\ldots, r_k \ge 0\;,\;k \ge 1\,. 
\label{dkDna}\ee

\bigskip 

\bn
{\em Proof of Theorem 1} \ Lemma 2 and Chebyshev's inequality imply that 
\be 
X_n(r) := \D_{B,n}^{(\a)}(r) - 2\,\frac{\a_{red}^{(2)}(c_n^{\a}\,r\,B)}{|W_n|^{\a}}\;\frac{N(W_n) - \l\,|W_n|}{\sqrt{|W_n|}} \Plongn 0\quad\mbox{for all}\;\; r \ge 0\,.
\label{Xn}\ee
In other words, $X_n(r)$ converges in probability to zero, that is, $\pr( |X_n(r)| \ge \e) \longn 0$ for all $\e > 0\,$. Now, we rearrange the terms in (\ref{Xn}) so that 
\bea  
\D_{B,n}^{(\a)}(r) &=& X_n(r) + Y_n(r) + 2\,\l\,|r\, B|\;\frac{N(W_n) -\l\,|W_n|}{\sqrt{|W_n|}}\,,\quad\qquad\qquad\qquad\label{DnXY}\\
\mbox{where}\qquad\qquad\qquad\qquad    Y_n(r) &:=& 2\,\g_{red}^{(2)}(c_n^{\a}\,r\,B)\;\frac{N(W_n) -\l\,|W_n|}{|W_n|^{\a+1/2}}\,.
\n\eea
Since $\a > 0$, it follows from  $\|\g_{red}^{(2)}\|_{var} < \infty$ and (\ref{cltN}) that $Y_n(r)$ also converges in probability to zero. 
Finally, condition (\ref{cltN}) and  Slutsky's lemma, see Billingsley (1968), Chapt. 1, provide the assertion (\ref{dDna}) of Theorem 1.
To prove (\ref{dkDna}) we use the Cram\'{e}r-Wold device, see Billingsley (1968), Chapt. 1, and show that, for any $\xi_1,\ldots,\xi_k \in \R^1$  
satisfying $\xi_1^2+\cdots+\xi_k^2 \ne 0$, the linear combination

\be 
\sum_{i=1}^k \xi_i\,\D_{B,n}^{(\a)}(r_i) = \sum_{i=1}^k \xi_i\,X_n(r_i) + \sum_{i=1}^k \xi_i\,Y_n(r_i) + \,2\,\l\,|B|\,
\frac{N(W_n)-\l\,|W_n|}{\sqrt{|W_n|}}\,\sum_{i=1}^k \xi_i\,r_i^d
\n\ee

\mn
converges in distribution to the Gaussian random variable $2\,\l\,|B|\,\sum_{i=1}^k \xi_i\,r_i^d\;{\mathcal N}(0,\sigma^2)\,$. As before in the
proof of (\ref{dDna}), this follows immediately by applying (\ref{cltN}), $X_n(r_i) \Plongn 0$ and $Y_n(r_i) \Plongn 0$ for $i = 1,\ldots, k$ 
together with Slutsky's lemma. $\Box$

\bn
\textbf{Corollary 1} \ The multivariate CLT (\ref{dkDna}) and (\ref{cltN}) give rise to define a $\chi^2$-type test statistic 

\be
T_{B,n}^{(\a)}(r_1,\ldots,r_k) := \frac{1}{\l^2\,\sigma^2} \sum_{i=1}^k \bigg(\frac{\D_{B,n}^{(\a)}(r_i)-\D_{B,n}^{(\a)}(r_{i-1})}{r_i^d-r_{i-1}^d}\bigg)^2 
+ \frac{|W_n|}{\sigma^2}\,\big(\,\widehat{\l}_n - \l\,\big)^2 \,.
\n\ee

\mn
For any $0 := r_0  < r_1 < \cdots < r_k < \infty\,$, $k \ge 1$ and all $\a \in (0,1)$ it holds 

\be 
T_{B,n}^{(\a)}(r_1,\ldots,r_k) \dlongn ( 4\,k\,|B|^2 + 1 )\;{\mathcal N}(0,1)^2\,, 
\n\ee

\mn
where the limit remains the same when $\l^2$ and $\sigma^2$ are replaced by their mean-square consistent estimators (\ref{l2est}) and (\ref{sigest}), respectively,
provided the assumptions of Lemma 4 are satisfied. 

\bn
{\em Proof of Corollary 1} \ Using (\ref{DnXY}) for $r \in \{r_1,\ldots,r_k\}$ gives

\be 
\frac{\D_{B,n}^{(\a)}(r_i)-\D_{B,n}^{(\a)}(r_{i-1})}{r_i^d-r_{i-1}^d} = 2\,\l\,|B|\,\sqrt{|W_n|}\,\big( \widehat{\l}_n - \l \big)
+ \frac{X_n(r_i) + Y_n(r_i)- X_n(r_{i-1}) - Y_n(r_{i-1})}{r_i^d-r_{i-1}^d}
\n\ee

\mn
for $i=1,\ldots,k$. By squaring the previous equality and summing up all squares  we get together with (\ref{cltN}) and  Slutsky's lemma the assertion 
of Corollary 1. $\Box$

\bn
\textbf{Remark 2} \ The CLT (\ref{cltN}) for the number $N(W_n)$ can be verified for various classes of stationary PPes, see Ivanoff (1982), among them Poisson cluster processes, see Heinrich (1988), or $\a$-determinantal PPes, see e.g. Heinrich (2016) and references therein. Note that for any Brillinger-mixing PP  condition (\ref{cltN}) is satisfied, see Jolivet (1981), Heinrich and Schmidt (1985), Karr (1987) or Biscio and Lavancier (2017).

\bigskip
It is noteworthy that the squared intensity $\l^2$ occurring in $\D_{B,n}^{(\a)}(r)$ cannot
be replaced by the asymptotically unbiased estimator (\ref{l2est}) without changing the limit as reveals  

\bn
\textbf{Remark 3} \ Under the conditions of Lemma 2 we have

\be
\lim_{n \to \infty}\,|W_n|^{1 - 2\a}\,\ew \bigl(\,(\widehat{\l^2 K_B})_n(c_n^{\a}\,r)\, -
\,(\widehat{\l^2})_n\,K_B(c_n^{\a}\,r)\,\bigr)^2\, = \, 0 \quad\mbox{for all}\quad r \ge 0\,.
\label{hatl2}\ee

\mn
The proof of (\ref{hatl2}) is based on the same arguments and techniques as used to prove Lemma 2. 
For the interested reader a detailed proof is given Appendix B.


\medskip
On the other hand, (\ref{cltN}) implies that $(\widehat{\l^2})_n$ is asymptotically normally distributed
which follows easily by applying  Slutsky's lemma:

\be 
\sqrt{|W_n|}\,\bigl(\,(\widehat{\l^2})_n - \l^2 \,\bigr) = \sqrt{|W_n|}\,\big(\,\widehat{\l}_n -\l\,\big)\,
\bigl(\,\widehat{\l}_n +\l\,\bigr) - \frac{\widehat{\l}_n}{\sqrt{|W_n|}}\; \dlongn \; 2\,\l\;{\mathcal N}(0,\sigma^2)\,.
\n\ee

\bn
\textbf{Theorem 2} \ In addition to the assumptions of Theorem 1, suppose
that (\ref{density}) is satisfied. Then, for $1/4 < \a < 1\,$, the weak convergence 

\be
\{\D_{B,n}^{(\a)}(r): r \in [0,R]\} \dlongn \{ G_B(r):= 2\,\l\,|B|\,r^d\;{\mathcal N}(0,\sigma^2) : r \in [0,R] \}
\label{Thm2}\ee

\mn
holds in the Skorohod-space $D[0,R]\,$.  If, in addition to (\ref{density}), the Lebesgue densities of $\g_{red}^{(3)}$ and $\g_{red}^{(4)}$
exist and satisfy (\ref{density34}), then  the weak convergence (\ref{Thm2}) in $D[0,R]$  holds for $0 < \a < 1\,$.

\bn
{\em Proof of Theorem 2} \ The multivariate CLT (\ref{dkDna}) states that  all finite-dimensional distributions of $(\D_{B,n}^{(\a)}(r))_{r \ge 0}$ converge to the corresponding distributions of the continuous mean zero Gaussian process  $(G_B(r))_{r \ge 0}$  with covariance function $\ew G_B(s)\,G_B(t) = 4\,\l^2\,\sigma^2\,|B|^2\,s^d\,t^d$. 
To verify the tightness of the sequence of random functions $\D_{B,n}^{(\a)} := (\D_{B,n}^{(\a)}(r))_{r \in [0,R]}$ in the Skorohod-space $D[0,R]\,$, see Billingsley (1968), Chapt. 3, we shall show that for each $\e > 0$ and $\eta > 0$ there exists a $\delta \in (0,1)$, such that

\be 
\pr(\,w(\D_{B,n}^{(\a)},\delta) \ge \e\,) \,\le \,\eta\quad\mbox{for all sufficiently large}\;n\,,
\label{w}\ee

\mn
where $w(f,\delta) :=\sup\{|f(t)-f(s)| : |t - s| \le \delta, 0 \le s,t \le R\}$ denotes the modulus of continuity of a so-called c\`{a}dl\`{a}g-functions 
$f:[0,R] \mapsto \R^1$ being right-continuous on $[0,R)$ and having limits from the left on $(0,R]$, i.e. $f \in D[0,R]\,$.
 
Since $r B$ is compact for all $r\ge 0$ and $s\,B \subseteq t\,B$ for $s \le t$,  the function $r \mapsto (\widehat{\l^2 K_B})_n(c_n^{\a}\,r)$
is piecewise constant, non-decreasing, right-continuous on $[0,R]$ ($\pr$-a.s.) with upward-jumps at $c_n^{-\a}\,\|X_i-X_j\|_B$ for $i \ne j$. Furthermore,  
$\l\,K_B(c_n^{\a}\,r) = \a_{red}^{(2)}(c_n^{\a}\,r\,B)$ is also non-decreasing and right-continuous (even absolutely continuous due (\ref{density})) 
for $r \ge 0\,$. Hence, by (\ref{KBest}) we have $\D_{B,n}^{(\a)}(0) = 0$ and $\pr(\D_{B,n}^{(\a)} \in D[0,R]) = 1$ for all $n \ge 1\,$. 

\medskip
Since $|(t B) \setminus (s B)| =(t^d - s^d)\,|B| \le d\,R^{d-1}\,|B|\,(t-s)$ for $0 \le s \le t \le R$,  Lemma 3 yields

\be 
\ew \bigl(\,\D_{B,n}^{(\a)}(t) - \D_{B,n}^{(\a)}(s)\,\bigr)^2\, \le \,a_1\,(d\,R^{d-1}\,|B|)^2\,(t - s)^2
+ d\,R^{d-1}\,|B|\,\frac{a_2\,(t - s)}{|W_n|^{\a}} + \frac{a_3}{|W_n|^{2\a}}\,.
\n\ee

\mn
For given $\e \in (0,1)$, for some $\beta \in (0,1)$ and sufficiently large $n_0=n_0(\e)$, we choose $s,t \in [0,R]$ such that 
$\e/|W_n|^{2\,\a} \le |t-s|^{1+\beta} \le 1$ for all $n \ge n_0$. Such a choice of $s,t \in [0,R]$ implies  
$|t-s|/|W_n|^{\a} \le |t-s|^{1+\beta}/\e$ and $(t-s)^2 \le |t-s|^{1+\beta}/\e$ so that 

\be 
\ew \bigl(\,\D_{B,n}^{(\a)}(t) - \D_{B,n}^{(\a)}(s)\,\bigr)^2\, \le a\;|t-s|^{1+\beta}/\e\;,
\label{inc1}\ee

\mn
where $a = a_1\,(d\,R^{d-1}\,|B|)^2 + a_2\,d\,R^{d-1} + a_3$ depends on $d, R, B, \l$ and $\|\g_{red}^{(k)}\|_{var}\,$, $k=2,3,4\,$.

Since $(\widehat{\l^2 K_B})_n(c_n^{\a}\,r)$ as well as $\l^2 K_B(c_n^{\a}\,r)$ is non-decreasing in $r \ge 0$, we find after a short 
calculation that, for any $s \le t \le u\,$,

\be 
|\D_{B,n}^{(\a)}(t) - \D_{B,n}^{(\a)}(s)| \le |\D_{B,n}^{(\a)}(u) - \D_{B,n}^{(\a)}(s)| + |W_n|^{1/2-\a}\,\l^2\,(K_B(c_n^{\a}\,u) - K_B(c_n^{\a}\,s))\,,
\label{inc2}\ee

\mn
where the second term on the r.h.s. is bounded by $a_0\,\l\,|W_n|^{1/2}\,|(u B) \setminus (s B)|$ due to (\ref{density}). 

\medskip
In the next step we determine an upper bound of $\sup_{s \le t \le s + m\,h}|\D_{B,n}^{(\a)}(t) - \D_{B,n}^{(\a)}(s)|$, 
where $m$ is a positive integer and the number $h$ is chosen such that $\e/|W_n|^{2\,\a} \le h^{1+\beta}$ and $h \le \min\{1, R\}\,$.
We consider the random variables $\zeta_k := \D_{B,n}^{(\a)}(s + k\,h)-\D_{B,n}^{(\a)}(s + (k-1)\,h)$ for $k=1,\ldots,m\,$ (for $s \le R$ and with $s + k\,h$ 
replaced by $R$ if it exceeds $R$). From (\ref{inc1}) it follows

\be 
\pr(|\D_{B,n}^{(\a)}(s + j\,h))-\D_{B,n}^{(\a)}(s + i\,h) | \ge \tau) = \pr(|\zeta_{i+1} + \cdots + \zeta_j| \ge \tau) \le \frac{a}{\e\,\tau^2}\,\big((j-i)\,h\big)^{1+\beta}   
\n\ee

\mn
for any $\tau > 0$ and $1 \le i \le j \le m\,$, where  $m$ does not exceed $\lfloor R/h \rfloor + 1\,$.  Now, we are in a position to apply Theorem 12.2 in Billingsley (1968), p. 94, which admits the estimate 

\be 
\pr(\max_{1 \le k \le m} |\D_{B,n}^{(\a)}(s + k\,h))-\D_{B,n}^{(\a)}(s)| \ge \tau) = \pr(\max_{1 \le k \le m}|\zeta_1 + \cdots + \zeta_k| \ge \tau)\,\le \,\frac{a_4\,a}{\e\,\tau^2}\,\big(m\,h\big)^{1+\beta}\,,  
\n\ee

\mn
where $a_4$ (occurring in Theorem 12.2) only depends on $\beta$.  Using (\ref{inc2}) for $u = t+h$ yields 

\be 
|\D_{B,n}^{(\a)}(t) - \D_{B,n}^{(\a)}(s)| \;\le\; |\D_{B,n}^{(\a)}(s+h) - \D_{B,n}^{(\a)}(s)| + a_5\,|W_n|^{1/2}\,h\quad\mbox{with}\;\;a_5 =  a_0\,\l\,|B|\,d\,R^{d-1}\,,
\n\ee

\mn
whence we conclude that

\be 
\sup_{s\le t \le s+m\,h}|\D_{B,n}^{(\a)}(t) - \D_{B,n}^{(\a)}(s)| \;\le\; 3\, \max_{1 \le k \le m} |\D_{B,n}^{(\a)}(s + k\,h))-\D_{B,n}^{(\a)}(s)|+ a_5\,|W_n|^{1/2}\,h\,.
\n\ee

To obtain the latter estimate we simply use the fact the sup is taken in at least one of the intervals $[s+(k-1)\,h, s+k\,h]$ for $k=1,\ldots,m\,$.
We specify now the choice of $h, m$ and $\delta \in (0,\min\{1,R\})$ for given (small enough) $\e > 0,\, \eta > 0$ as follows:
Since $\a > 1/4$ we may put $\beta = 4\,\a - 1 > 0$ for $\a < 1/2$ (and $\beta = \a$ for $\a \in [1/2,1)$) and hence $h$ can be chosen to fulfill the inequality 
$\e^{\frac{1}{1+\beta}}/|W_n|^{\frac{2\,\a}{1+\beta}} \le h < \e/a_5\,|W_n|^{1/2}$  for $n \ge n_0(\e)\,$. This choice of $h$ combined with the previous estimates for 
 $\tau =\e$  gives 

\be 
\pr(\sup_{s\le t \le s+m\,h}|\D_{B,n}^{(\a)}(t) - \D_{B,n}^{(\a)}(s)|\ge 4\,\e)\, \le\, \pr(\max_{1 \le k \le m} |\zeta_1 + \cdots + \zeta_k|\ge \e) 
\,\le \,\frac{a_4\,a}{\e^3}\,\big(m\,h\big)^{1+\beta} \,. 
\n\ee

\medskip
Setting $m = \lfloor (\e^3\,\eta)^{1/\beta}/h \rfloor$ and  $\delta  = m\,h$ (such that $0 < \delta \le (\e^3\,\eta)^{1/\beta} < \min\{1,R\}$) we get that

\be 
\pr(\sup_{s\le t \le s+\delta}|\D_{B,n}^{(\a)}(t) - \D_{B,n}^{(\a)}(s)|\ge 4\,\e)\,\le\,\delta\;\frac{a_4\,a}{\e^3}\,\big(m\,h\big)^{\beta} \le a_4\,a\,\delta\;\eta
\quad\mbox{for}\;\; n \ge n_0(\e)\;,\; s \in [0, R - \delta]\,. 
\n\ee

Finally, applying the Corollary to Theorem 8.3 in Billingsley (1968) we obtain that for any  $\e, \eta > 0$ there exists a $\delta=\delta(\e,\eta) > 0$ such that   

\be 
\pr(w(\D_{B,n}^{(\a)},\delta)\ge 12\,\e) \le (1 + \lfloor R/\delta \rfloor)\, a_4\,a\,\delta\,\eta \le  a_4\,a\,(1+R)\,\eta \quad\mbox{for all sufficiently large}\;\;n\,,
\n\ee

\mn
which coincides with (\ref{w}) up to the factors 12 and $a_4\,a\,(1+R)$, respectively. Therefore, the sequence $(\D_{B,n}^{(\a)})_{n\ge 1}$ is tight in
$D(0,R)$ and together with the convergence of its the finite-dimensional distributions the proof of (\ref{Thm2}) for $1/4 < \a < 1$ is complete. The validity of 
(\ref{Thm2}) can be extended to $0 < \a \le 1/4$ under the additional assumption (\ref{density34}). To see this, let $\a \in (0,1)$. For given $\e \in (0,1)$ 
 we choose $s,t \in [0,R]$ such that $\e/|W_n|^{\a} \le |t-s|^{\a} \le 1$ for all $n \ge n_0$.  In view of (\ref{rem1}) such a choice of $s,t \in [0,R]$
 enables us to replace  (\ref{inc1}) by 

\be 
\ew \bigl(\,\D_{B,n}^{(\a)}(t) - \D_{B,n}^{(\a)}(s)\,\bigr)^2\, \le b\;|t-s|^{1+\a}/\e\quad\mbox{with}\quad b = a_1\,(d\,R^{d-1}\,|B|)^2 + a_2\,d\,R^{d-1} \;,
\label{inc3}\ee

\mn
where the constants $a_1$, $a_2$ are from (\ref{rem1}). Now, the increment $h$ can be chosen to fulfill the inequality 
$\e^{\frac{1}{\a}}/|W_n| \le h < \e/a_5\,|W_n|^{1/2}$  for $n \ge n_0(\e)\,$. The remaining steps are the same as before with $\beta = \a$. 
This terminates the proof of Theorem 2. $\Box$

\mn

\bn
\textbf{Corollary 2} \ Let $V | [0,R] \mapsto [0,\infty)$ be non-decreasing, right-continuous  function such that $V(R) < \infty$. 
Under the assumptions of Theorem 2 (additionally with (\ref{density34})) we have for $1/4 < \a \le 1/2\;$ ($\,0 < \a \le 1/2\,$ )

\be 
\int_0^R\,\bigl(\,\D_{B,n}^{(\a)}(r)\,\bigr)^2\,{\rm d}V(r) 
\dlongn \;4\,\l^2\,|B|^2\;\int_0^R\,r^{2d}\,{\rm d}V(r)\;{\mathcal N}(0,\sigma^2)^2\,.  
\label{cor2}\ee

\bn
\textbf{Corollary 3} \ Let $v$ be a c\`{a}dl\`{a}g-function on $[0,R]$, i.e. $v \in D[0,R]$. Under the assumptions  of Theorem 2 
(additionally with (\ref{density34})) we have for $1/4 < \a \le 1/2\;$ ($\,0 < \a \le 1/2\,$ )

\be
\sup_{0 \le r \le R}\,\big|\,\D_{B,n}^{(\a)}(r)\,v(r)\,\big|
\dlongn\; 2\,\l\,|B|\;\sup_{0 \le r \le R}|r^d\,v(r)|\;|\,{\mathcal N}(0,\sigma^2)\,|\,. 
\label{cor3}\ee

\bn
The proofs of Corollary 2 and 3 rely on the fact that both functionals $F_1(f) = \int_0^R f^2(r)\,{\rm d}V(r)$ and  $F_2(f) = \sup_{0 \le r \le R}|f(r)\,v(r)|$ 
are continuous on $D[0,R]$ in the following sense: Let $(f_n)$ be a sequence in $D[0,R]$ having a continuous limit $f$ w.r.t. to the Skorohod metric in $D[0,R]$,
see Chapt. 3 in Billingsley (1968) for a precise definition.  By virtue of this definition it can be shown that $F_i(f_n) \longn F_i(f)$ for $i=1,2$. Since the 
limit process $\{G_B(r) : r \in [0,R]\}$ in Theorem 2 is ($\pr$-a.s.) continuous, a direct application of the {\em continuous mapping theorem}, see Theorem 5.1 
in Billingsley (1968), yields the limit theorems (\ref{cor2}) and (\ref{cor3}). $\Box$  

\bn
\textbf{Remark 4} \ The trivial weight functions  $V(r) = r$ in (\ref{cor2}) and  $v(r) = 1$ in (\ref{cor3}) provide analogues
to the classical Cram\'{e}r-von Mises  and Kolmogorov-Smirnov test statistic. By means of the functions $V$ and $v$ one can place
higher weights on $|\,\D_{B,n}^{(\a)}(r)\,|$ in particular for small $r$-values, e.g. $v(r) = e^{-a\,r}$ for $a \ge d/R$ or 
$V'(r) = e^{-b\,r^{2d+1}}$ for $b > 0$.    

\bn
\section{Goodness-of-fit-tests for checking point process hypotheses}

\bn
{\bf Aim:} Testing the null hypothesis $H_0 : P = P_0$ vs. $H_1 : P \ne P_0$ by checking the goodness-of-fit 
of the (known) $K-$function $K_{0,B}(\cdot)$ and intensity $\l_0$ of $N \sim P_0$.

\bn
For this purpose we change the notation (\ref{Dna}) and put   

\be 
\D^{(\a)}_{0,B,n}(r) := |W_n|^{1/2 - \a}\,\big(\,(\widehat{\l^2 K_B})_n(c_n^{\a}\,r) - \l_0^2\,K^{}_{0,B}(c_n^{\a}\,r)\,\big)\quad\mbox{for}\quad 0\le r \le R
\;\;,\;\;0 \le \a < 1\,.
\n\ee 

\bn
As consequence of (\ref{remainder})  in order to distinguish between PPes with different $K_B-$function, the below Theorems 3 - 5 are 
formulated only  for $0 < \a \le 1/2\,$. In the particular case $\a = 1/2$ the normalizing sequence $|W_n|^{1/2 - \a}$ is constant equal to 1. 

First we establish a $\chi^2$-type test as an immediate outcome of Corollary 1. According to this under the null hypothesis $H_0$ and provided that
$N \sim P_0$ satisfies the conditions of Theorem 1 and Lemma 4, the test statistic

\be
\widehat{T_{0,B,n}^{(\a)}}(r_1,\ldots,r_k) :=  \frac{1}{(\widehat{\l^2})_n\,(\widehat{\sigma^2})_n} \sum_{i=1}^k 
\bigg(\frac{\D_{0,B,n}^{(\a)}(r_i)-\D_{0,B,n}^{(\a)}(r_{i-1})}{r_i^d-r_{i-1}^d}\bigg)^2  + \frac{|W_n|\,(\widehat{\l}_n - \l_0)^2}{(\widehat{\sigma^2})_n} 
\n\ee

\mn
possesses the limit distribution function $F_{k,B}(x) := \pr((4\,k\,|B|^2 + 1)\,{\mathcal N}(0,1)^2 \le x)$ (as $n \to \infty$) 
for any $0 := r_0  < r_1 < \cdots < r_k < \infty\,$, $k \ge 1$ and all $\a \in (0,1/2]\,$. Given a {\em significance level} $\g > 0$ ( $=$ probability 
of the {\em type I error} ), $H_0$ is accepted (for sufficiently large $n$) if  $\widehat{T_{0,B,n}^{(\a)}}(r_1,\ldots,r_k) \le (4\,k\,|B|^2 + 1)\,z_{1-\g/2}^2\,$, where 
$z_p$ denotes the $p$-quantile of the standard normal distribution $\Phi(x) = \pr({\mathcal N}(0,1) \le x)$. This follows from   
 $F_{k,B}((4\,k\,|B|^2 + 1)\,z_{1-\g/2}^2) = 1 - \pr(|{\mathcal N}(0,1)| > z_{1-\g/2} )= 1 - \g\,$.  

\mn

\bn
\textbf{Theorem 3} \ Let $N_0 \sim P_0$ satisfy the conditions of Theorem 2 (additionally with (\ref{density34})) and Lemma 4.
Then, under the null hypothesis $H_0\,$,  we have for $1/4 < \a \le 1/2\;$ ($\,0 < \a \le 1/2\,$ )

\be 
\frac{\int_0^R\,\big(\,\D^{(\a)}_{0,B,n}(r)\,\big)^2\,{\mathrm d}V(r)}{(\widehat {\l^2})_n\,
(\widehat {\sigma^2})_n} + \frac{|W_n|\,(\,\widehat{\l}_n - \l_0\,)^2}{(\widehat{\sigma^2})_n}
\dlongn \Big(\,4\,|B|^2\;\int_0^R\,r^{2d}\,{\rm d}V(r) + 1\,\Big)\;{\mathcal N}(0,1)^2
\label{T3CM}\ee

and
\be 
\frac{\sup\limits_{0 \le r \le R}\big|\,\D^{(\a)}_{0,B,n}(r)\,v(r)\,\big|}{\widehat{\l}_n\,\sqrt{(\widehat {\sigma^2})_n}}
+ \frac{\sqrt{|W_n|}\,|\,\widehat {\l}_n - \l_0\,|}{\sqrt{(\widehat{\sigma^2})_n}}
\dlongn \Big(\,2\,|B|\,\sup_{0 \le r \le R}|r^d\,v(r)| + 1\,\Big)\;|\,{\mathcal N}(0,1)\,|\,,
\label{T3KS}\ee

\mn
where the weight functions $V$ and $v$ are the same as in Corollary 2 and 3, respectively.

\bn
{\em Proof of Theorem 3} \ In view of Slutsky's lemma it suffices to prove (\ref{T3CM}) when  $(\widehat {\l^2})_n$ and 
$(\widehat {\sigma^2})_n$ are replaced  by $\l_0^2$ and  $\sigma_0^2$ ( = asymptotic variance of $N(W_n)/\sqrt{|W_n|}$ 
under $H_0$).
   
Using the relations (\ref{Xn}) and (\ref{DnXY}) under $H_0 : N \sim P_0$ we obtain after some rearrangements the following identity

\bea 
\frac{1}{\l_0^2}\int_0^R\big( \D^{(\a)}_{0,B,n}(r) \big)^2\,{\mathrm d}V(r) +  |W_n|\,(\widehat{\l}_n - \l_0)^2 = 
\Big(\frac{N(W_n) -\l_0\,|W_n|}{\sqrt{|W_n|}}\Big)^2\,\Big( \int_0^R 4\,|r\, B|^2\,{\mathrm d}V(r) + 1 \Big) &&\n\\ 
&&\n\\
+ \int_0^R\Big(\frac{X_n(r) + Y_n(r)}{\l_0}\Big)\,\Big(\frac{X_n(r) + Y_n(r)}{\l_0} + 
4\,|r\, B|\;\frac{N(W_n) -\l_0\,|W_n|}{\sqrt{|W_n|}}\Big)\,{\mathrm d}V(r)\,.\phantom{00000000}&&
\label{intCM}\eea

Now, we apply the tightness of the sequence $\{ \D^{(\a)}_{0,B,n}(r) : r\in [0,R] \}$ which implies that  
 $\sup_{0 \le r \le R}|X_n(r)| \Plongn 0$ and  $\sup_{0 \le r \le R}|Y_n(r)| \Plongn 0$. Thus, the integral in (\ref{intCM}) tends 
 to zero in probability.  Finally, the assumption (\ref{cltN}) combined once more with Slutsky's lemma  shows  the normal convergence 
 
\be  
\frac{1}{\l_0^2}\int_0^R\big( \D^{(\a)}_{0,B,n}(r) \big)^2\,{\mathrm d}V(r) + |W_n|\,(\widehat{\l}_n - \l_0)^2 \dlongn  
\Big( \int_0^R 4\,|r\, B|^2\,{\mathrm d}V(r) + 1 \Big)\;{\mathcal N}(0,\sigma_0^2)^2
\n\ee

\mn 
which is equivalent to (\ref{T3CM}). The proof of the second assertion (\ref{T3KS}) relies on the same  arguments as before.
The relations  (\ref{Xn}) and (\ref{DnXY}) (under  $H_0 : N \sim P_0$) allows to derive the ($\omega$-wise) estimate 

\bea 
\Bigg| \frac{\sup\limits_{0 \le r \le R}\big|\,\D^{(\a)}_{0,B,n}(r)\,v(r)\,\big|}{\l_0} + \sqrt{|W_n|}\,\big|\,\widehat {\l}_n - \l_0\,\big| 
- \big(\,2\,|B|\,\sup_{0 \le r \le R} |r^d\,v(r)| + 1\,\big)\,\frac{\big|\,N(W_n) -\l_0\,|W_n|\,\big|}{\sqrt{|W_n|}} \Bigg| &&\n\\ 
\phantom{0000}\n\\
\le  \frac{1}{\l_0}\,\sup\limits_{0 \le r \le R}\big|\,(X_n(r) + Y_n(r))\,v(r)\,\big|\; \Plongn 0 \,.\phantom{0000000000000000000}&&
\n\eea

Together with (\ref{cltN}) and Slutsky's lemma  it follows the distributional  convergence 
\be
\frac{\sup\limits_{0 \le r \le R}\big|\,\D^{(\a)}_{0,B,n}(r)\,v(r)\,\big|}{\l_0} + \sqrt{|W_n|}\,\big|\,\widehat {\l}_n - \l_0\,\big| 
\dlongn \big(\,2\,|B|\,\sup_{0 \le r \le R} |r^d\,v(r)| + 1\,\big)\,|\,{\mathcal N}(0,\sigma_0^2)\,|\,,
\n\ee
which turns out to be equivalent to (\ref{T3KS})). Thus, the proof of Theorem 3 is complete. $\Box$

\bigskip
A further application of the Theorem 2 consists in testing the hypothesis that two independent PPes $N_a \sim P_a$
and $N_b \sim P_b$ possess the same distribution, that is, testing $H_0^{a,b} : P_a = P_b$ vs. $H_1^{a,b} : P_a \ne P_b$.
Note that the parameters $\l_a$, $\sigma_a^2$ and the $K$-function $K_{a,B}$ and their estimators 
are related to $P_a$ and an observation of $N_a$ in $W_n$  which also applies for quantities with subscript $b$.
It should be mentioned that first attempts to tackle this two-sample test problem for stationary PPes were made in Doss (1989).

\bn
\textbf{Theorem 4} \ Let  $N_a \sim P_a$ and $N_b \sim P_b$ be two independent $B_4$-mixing PPes
on $\R^d$ which can be observed in $W_n$ and satisfy the conditions of  Theorem 2 (additionally with (\ref{density34})) and Lemma 4. 
Then, under the null hypothesis $H_0^{a,b}\,$,  we have for $1/4 < \a \le 1/2\;$ ($\,0 < \a \le 1/2\,$ )

\bea  
&&\frac{|W_n|\,\int_0^R\,\bigl(\,(\widehat {\l_a^2 K_a^{}})^{}_n(c_n^{\a}\, r)\, - \,
(\widehat{\l_b^2 K_b^{}})^{}_n(c_n^{\a}\, r)\,\bigr)^2\,{\rm d}V(r)}{|W_n|^{2\,\a}\,\big((\widehat {\l_a^2})_n^{}\,
(\widehat{\sigma_a^2})_n^{}\, + \,(\widehat {\l_b^2})_n^{}\,(\widehat {\sigma_b^2})_n^{}\big)} + 
\frac{|W_n|\,\big(\,(\widehat {\l_a})_n - (\widehat {\l_b})_n\,\big)^2}{(\widehat {\sigma_a^2})^{}_n\,+\,(\widehat {\sigma_b^2})^{}_n\,}\phantom{00} \n\\
&&\label{T4CM}\\
&&\n\\
&&\phantom{00000000000}\dlongn \;\Big(\,4\,|B|^2\,\int_0^R\,r^{2d}\,{\rm d}V(r) + 1\,\Big)\;{\mathcal N}(0,1)^2\,.
\n\eea

\mn
{\em Proof of Theorem 4} \ In case of the PP  $N_a \sim P_a$ we use the abbreviation
\be 
\D^{(\a)}_{a,B,n}(r) := |W_n|^{1/2 - \a}\,\big(\,(\widehat{\l_a^2 K_{a,B}})_n(c_n^{\a}\,r) - \l_a^2\,K^{}_{a,B}(c_n^{\a}\,r)\,\big)\quad\mbox{for}\quad 0 \le r \le R
\n\ee 
and correspondingly $\D^{(\a)}_{b,B,n}(r)$ for  $N_b \sim P_b\,$. Further, define $Q_{a,n} := (N_a(W_n)-\l_a\,|W_n|)/\sqrt{|W_n|}$ and likewise $Q_{b,n}$.  
Using the relations (\ref{Xn}) and (\ref{DnXY}) for both PPes $N_a$ and $N_b$ (with obvious changes of the notation) we arrive at

\be 
\D^{(\a)}_{a,B,n}(r)- \D^{(\a)}_{b,B,n}(r) = Z_n^{a,b}(r) + 2\,| r\,B |\,(\,\l_a\,Q_{a,n} - \l_b\,Q_{b,n}\,)\,,
\label{Zabn}\ee

\mn
where $Z_n^{a,b}(r) :=  X_{a,n}(r) - X_{b,n}(r) + Y_{a,n}(r) - Y_{b,n}(r)\,$. Under $H_0^{a,b}$ (which implies $\l_a = \l_b$ and $K_{a,B} = K_{b,B}$) we may use the ($\omega$-wise) equality $\sqrt{|W_n|}\,\big(\,(\widehat {\l_a})_n - (\widehat {\l_b})_n\,\big) = Q_{a,n} - Q_{b,n}\,$.
After squaring, integrating w.r.t. ${\rm d}V$ and dividing both sides of (\ref{Zabn}) by $\l_a^2\,\sigma_a^2 + \l_b^2\,\sigma_b^2$ and provided  $H_0^{a,b} : P_a = P_b$ is true, we obtain  the equality 

\bea
\frac{|W_n|\,((\widehat {\l_a})_n - (\widehat {\l_b})_n)^2}{\sigma_a^2 + \sigma_b^2} + \int_0^R \frac{(\D^{(\a)}_{a,B,n}(r) - \D^{(\a)}_{b,B,n}(r))^2}{\l_a^2\,\sigma_a^2 + \l_b^2\,\sigma_b^2}\,{\rm d}V(r)  =  \frac{(Q_{a,n} - Q_{b,n})^2}{\sigma_a^2 + \sigma_b^2}\phantom{0000000} && \n\\
&&\n\\
&&\n\\
+ 4\,\int_0^R |r B|^2\,{\rm d}V(r) \,\frac{(Q_{a,n} - Q_{b,n})^2}{\sigma_a^2 + \sigma_b^2} + \int_0^R \frac{Z_n^{a,b}(r)\,\big(\,Z_n^{a,b}(r) + 4\,|r B|\,(\l_a\,Q_{a,n} - \l_b\,Q_{b,n})\,\big)}{\l_a^2\,\sigma_a^2 + \l_b^2\,\sigma_b^2}\,{\rm d}V(r)\,, &&
\n\eea

\mn
where  $\D^{(\a)}_{a,B,n}(r) - \D^{(\a)}_{b,B,n}(r)$ can be replaced by $|W_n|^{1-2\,\a}\,(\,(\widehat {\l_a^2 K_a^{}})^{}_n(c_n^{\a}\,r) - (\widehat{\l_b^2 K_b^{}})^{}_n(c_n^{\a}\, r)\,)$ for all $0 \le r \le R\,$. 
The assumptions of Theorem 2 put on both PPes $N_a \sim P_a$ and $N_b \sim P_b$ and their independence imply 

\be
\sup_{0 \le r \le R}|Z_n^{a,b}(r)| \Plongn 0\quad\mbox{and}\quad \frac{Q_{a,n} - Q_{b,n}}{\sqrt{\sigma_a^2 + \sigma_b^2}} \dlongn {\mathcal N}(0,1)\,, 
\n\ee
\mn
whence in summary we conclude that

\bea 
&&\frac{|W_n|^{1-2\,\a}}{\l_a^2\,\sigma_a^2 + \l_b^2\,\sigma_b^2}\,\int_0^R(\,(\widehat {\l_a^2 K_a^{}})^{}_n(c_n^{\a}\,r) - (\widehat{\l_b^2 K_b^{}})^{}_n(c_n^{\a}\,r)\,)^2\,{\rm d}V(r) + \frac{|W_n|\,((\widehat {\l_a})_n - (\widehat {\l_b})_n)^2}{\sigma_a^2 + \sigma_b^2}\n\\
&& \phantom{0000000000} \dlongn  \Big( 4\,|B|^2\,\int_0^R r^{2d}\,{\rm d}V(r) + 1 \Big)\;{\mathcal N}(0,1)^2\,.
\n\eea

\medskip 
After replacing the unknown parameters $\l_a^2, \l_b^2$ and $\sigma_a^2, \sigma_b^2$ by the mean-square consistent estimators given in (\ref{l2est}) and (\ref{sigest}), respectively, and  applying of Slutsky's lemma, the proof of Theorem 4 is complete. $\Box$

\bn
\textbf{Theorem 5} \ Let  $N_a \sim P_a$ and $N_b \sim P_b$ be two independent $B_4$-mixing PPes on $\R^d$ which 
can be observed in $W_n$ and satisfy the conditions of  Theorem 2 (additionally with (\ref{density34})) and Lemma 4. 
Then, under the null hypothesis $H_0^{a,b}\,$,  we have for $1/4 < \a \le 1/2\;$ ($\,0 < \a \le 1/2\,$ )

\bea 
&& \frac{\sqrt{|W_n|}\,\sup\limits_{0 \le r \le R}\,\big|\,( (\widehat{\l_a^2 K_a^{}})_n^{}(c_n^{\a}\,
r)\, - \,(\widehat{\l_b^2 K_b^{}})_n^{}(c_n^{\a}\, r) )\,v(r)\,\big|}{|W_n|^{\a}\,\sqrt{\,(\widehat {\l_a^2})^{}_n\,
(\widehat {\sigma_a^2})^{}_n \,+\, (\widehat {\l_b^2})^{}_n\,(\widehat {\sigma_b^2})^{}_n\,}} + 
\frac{\phantom{\int\limits_i}\sqrt{|W_n|}\,\big|\,(\widehat {\l_a})_n-(\widehat {\l_b})_n\,\big|}{\sqrt{\,
(\widehat {\sigma_a^2})^{}_n\,+\,(\widehat {\sigma_b^2})^{}_n\,}}\phantom{000} \n\\ 
&&\label{T5KS} \\
&&\n\\
&&\phantom{000000000000}\dlongn  \;\Big(\,2\,|B|\,\sup_{0 \le r \le R}|\,r^d\,v(r)\,| + 1\,\Big)\;|\,{\mathcal N}(0,1)\,|\,.
\n\eea

\mn
{\em Proof of Theorem 5} \ With the notation introduced  in the proof of Theorem 4 and (\ref{Zabn}) it is easily seen that

\be  
\Bigg| \frac{\sup\limits_{0 \le r \le R}\big|\,\big(\D^{(\a)}_{a,B,n}(r) - \D^{(\a)}_{b,B,n}(r)\big)\,v(r)\,\big|}{\sqrt{\l_a^2\,\sigma_a^2 + \l_b^2\,\sigma_b^2}} 
- 2\,|B|\,\sup_{0 \le r \le R} |r^d\,v(r)|\,\frac{\big|\,\l_a\,Q_{a,n} -\l_b\,Q_{b,n}\,\big|}{\sqrt{\l_a^2\,\sigma_a^2 + \l_b^2\,\sigma_b^2}} \Bigg| 
\n\ee
is bounded by the sequence $\sup\limits_{0 \le r \le R}\big|\,Z_n^{a,b}(r)\,v(r)\,\big|/\sqrt{\l_a^2\,\sigma_a^2 + \l_b^2\,\sigma_b^2}$ tending to zero in probability.  
Obviously, just as in the foregoing proof, if $H_0^{a,b}$ is true and if
(\ref{cltN}) is satisfied for both of the independent PPes $N_a$ and $N_b$ we have 

\be   
\frac{\sqrt{|W_n|}\,\big|(\widehat {\l_a})_n - (\widehat {\l_b})_n\big|}{\sqrt{\sigma_a^2 + \sigma_b^2}}
= \frac{\big|\,\l_a\,Q_{a,n} -\l_b\,Q_{b,n}\,\big|}{\sqrt{\l_a^2\,\sigma_a^2 + \l_b^2\,\sigma_b^2}} =
\frac{\big|\,Q_{a,n} - Q_{b,n}\,\big|}{\sqrt{\sigma_a^2 + \sigma_b^2}} \dlongn |\,{\mathcal N}(0,1)\,|\,.
\n\ee

Hence, if $H_0^{a,b}$ is true and all the assumptions of Theorem 4 are fulfilled,  

\bea 
&&\frac{\sqrt{|W_n|}\,\sup\limits_{0 \le r \le R}\big|\,( (\widehat{\l_a^2 K_a^{}})_n^{}(c_n^{\a}\,r)\, - 
\,(\widehat{\l_b^2 K_b^{}})_n^{}(c_n^{\a}\, r) )\,v(r)\,\big|}{|W_n|^{\a}\,\sqrt{\l_a^2\,\sigma_a^2 + \l_b^2\,\sigma_b^2}} 
+ \frac{\sqrt{|W_n|}\,\big|(\widehat {\l_a})_n - (\widehat {\l_b})_n\big|}{\sqrt{\sigma_a^2 + \sigma_b^2}}\n\\
&&\n\\
&& \phantom{0000000000000000000} \dlongn 
\Big( 2\,|B|\,\sup_{0 \le r \le R} |r^d\,v(r)| + 1 \Big)\, |\,{\mathcal N}(0,1)\,|\,.
\n\eea

\mn
The arguments at the of the proof of Theorem 4 confirm the assertion stated in Theorem 5. $\Box$

\bn
\textbf{Remark 5} \ The observation windows of $N_a$ and $N_b$ need not to be the same. Theorems 4 and 5 remain valid when the independent 
PPes $N_a$ and $N_b$ are given in $(W_n^a)$ and $(W_n^b)$ (each forming a CAS) such that $|W_n^a| = |W_n^b|$ for sufficiently large $n$. 

\bn
\textbf{Remark 6} \
Let there be given a significance level $\g > 0$ ( $=$ probability of the  type I error). Then, $H_0$ resp. $H_0^{a,b}$ is accepted (for large $n$) 
if  the test statistic on the l.h.s. of (\ref{T3CM}) resp. (\ref{T4CM}) does not exceed $\big( 4\,|B|^2\,\int_0^Rr^{2d}{\rm d}V + 1 \big)\,z_{1-\g/2}^2$ 
(or if the test statistic on the l.h.s. of (\ref{T3KS}) resp. (\ref{T5KS}) does not exceed $\big( 2\,|B|\,\sup_{r\in [0,R]}|r^d\,v(r)| + 1 \big)\,z_{1-\g/2}\,$), 
where $z_p$ is defined by $\Phi(z_p) = p\,$.

\bn
\section{Proofs of Lemma 2 and Lemma 3}

\mn
{\em Proof of Lemma 2} \ We start with the following decomposition

\bea 
&& \ew \Bigl[\,|W_n|\,(\widehat{\l^2 K_B})_n(c_n^{\a}\,r) - N(W_n)\,\a_{red}^{(2)}(c^{\a}_n\,r\,B)
\,-\,\a_{red}^{(2)}(c^{\a}_n\,r\,B)\,\bigl( N(W_n) - \l\,|W_n|\bigr)\,\Bigr]^2\n\\
&&\phantom{00}\n\\
&=& \bigl( \a_{red}^{(2)}(c^{\a}_n\,r\,B) \bigr)^2\,\var(N(W_n))
\,+\,{\ew}\Bigl[\,|W_n|\,(\widehat{\l^2 K_B})_n(c_n^{\a}\,r) - N(W_n)\,\a_{red}^{(2)}(c^{\a}_n\,r\,B)\,\Bigr]^2\n\\
&&\phantom{00}\n\\
&-& 2\,\a_{red}^{(2)}(c^{\a}_n\,r\,B)\,{\ew}\Bigl[\,\bigl( N(W_n) - \l\,|W_n|\bigr)\,
\Bigr[\,|W_n|\,(\widehat{\l^2 K_B})_n(c_n^{\a}\,r) - N(W_n)\,\a_{red}^{(2)}(c^{\a}_n\,r\,B)\,\Bigr]\n\\
&&\phantom{00}\n\\
&=& T^{(1)}_n(r\,B) + T^{(2)}_n(r\,B) - T^{(3)}_n(r\,B)\;.
\n\eea

\mn
The assertion of Lemma 2  follows from the three limits  

\be
\lim_{n \to \infty} |W_n|^{-(1+2\a)}\,T^{(1)}_n(r\,B) = \lim_{n \to \infty} |W_n|^{-(1+2\a)}\,T^{(2)}_n(r\,B) 
= \l^3\,|r\,B|^2\,\bigl(1+\g_{red}^{(2)}(\R^d)\bigr)\quad
\label{T12}\ee
and
\be
\lim_{n \to \infty} |W_n|^{-(1+2\a)}\,T^{(3)}_n(r\,B) = 2\,\l^3\,|r\,B|^2\,\bigl(1+\g_{red}^{(2)}(\R^d)\bigr)\,.
\label{T3}\ee

\medskip
Using the second-order stationarity of $N = \sum_{i\ge 1}\d_{X_i}$ and the definitions of $\g^{(2)}$ and
$\g_{red}^{(2)}$ we get      
\bea 
\var(N(W_n)) = {\ew}\bigl( N(W_n)\,-\,\l\,|W_n| \bigr)^2 &=& \l\,|W_n| + \g^{(2)}(W_n \times W_n) = \l\,|W_n|  + \n\\
&&\phantom{00}\n\\
\l \int_{\R^d}\int_{\R^d}\1_{W_n}(x)\,\1_{W_n}(y+x)\g_{red}^{(2)}({\rm d}y) {\rm d}x 
&=& \l\,|W_n| + \l\,\int_{\R^d}| W_n \cap (W_n-y) |\g_{red}^{(2)}({\rm d}y) \,,
\n\eea
whence in view of $\|\g_{red}^{(2)}\|_{var} < \infty$, Lemma 1 and Lebesgue's convergence theorem it follows that

\be
\frac{\a_{red}^{(2)}(c^{\a}_n\,r\,B)}{|W_n|^{\a}} = \l\,|r\,B| + \frac{\g_{red}^{(2)}(c^{\a}_n\,r\,B)}{|W_n|^{\a}}
\longn \;\l\,|r\,B| \quad\mbox{for}\quad \a > 0
\label{ared}\ee
and
\be
\frac{\g^{(2)}(W_n \times W_n)}{|W_n|} = \frac{\l}{|W_n|}\,\int_{\R^d} |W_n \cap (W_n - y)|\,\g_{red}^{(2)}({\mathrm d}y) 
\longn \;\l\,\g_{red}^{(2)}(\R^d)\,.
\n\ee
The latter two relations prove the first limit of (\ref{T12}).

\bigskip
Taking into account the unbiasedness of $(\widehat{\l^2 K_B})_n(c_n^{\a}\,r)$, i.e. $\ew (\widehat{\l^2 K_B})_n(r) 
= \l^2\,K_B(r) = \l\,\a_{red}^{(2)}(r\,B)\,$, and $\ew N(W_n) = \l\,|W_n|$ we can write $T^{(3)}_n(r\,B)$ as follows:
\bea 
T^{(3)}_n(r\,B) &=& 2\,\a_{red}^{(2)}(c^{\a}_n\,r\,B)\,\bigl[|W_n|\,\ew N(W_n)\,(\widehat{\l^2 K_B})_n(c_n^{\a}\,r) 
- \ew N^2(W_n)\,\a_{red}^{(2)}(c_n^{\a}\,r\,B)\,\bigr]\n\\
&&\phantom{00}\n\\
&=& 2\,\a_{red}^{(2)}(c^{\a}_n\,r\,B)\,|W_n|\,\bigl[\,\ew N(W_n)\,(\widehat{\l^2 K_B})_n(c_n^{\a}\,r)
- \l^2\,|W_n|\,\a_{red}^{(2)}(c_n^{\a}\,r\,B)\,\bigr] \n\\
&&\phantom{00}\n\\ 
&-& 2\,\big(\a_{red}^{(2)}(c^{\a}_n\,r\,B)\big)^2\,\var(N(W_n))\,.
\n\eea

\mn
Combining (\ref{ared}) and $\var(N(W_n))/|W_n| \longn \l\,\bigl(1 + \g_{red}^{(2)}(\R^d)\bigr)$ reveals that (\ref{T3}) 
is equivalent to

\be 
\lim_{n \to \infty}|W_n|^{-\a}\,\bigl[\,{\ew}N(W_n)\,(\widehat{\l^2 K_B})_n(c_n^{\a}\,r)- \l^2\,|W_n|\,
\a_{red}^{(2)}(c_n^{\a}\,r\,B)\,\bigr] = 2\,\l^2\,|r\,B|\,\bigr(1+\,\g_{red}^{(2)}(\R^d)\bigr)\,.
\label{rest}\ee

\noindent
For this purpose we rewrite $\ew N(W_n)\,(\widehat{\l^2 K_B})_n(c_n^{\a}\,r)$ as follows

\bea
&& \ew N(W_n)\,(\widehat{\l^2 K_B})_n(c_n^{\a}\,r) = \ew \sum_{i,j,k \ge 1 \atop i \ne j}\frac{\1_{W_n}(X_i)\,\1_{W_n}(X_j)\,
\1_{W_n}(X_k)}{|(W_n-X_i) \cap (W_n - X_j)|}\,\1_{c_n^{\a}\,r\,B}(X_j-X_i) \n\\
&&\phantom{00}\n\\
&=& \ew \Su_{i,j,k \ge 1} \frac{\1_{W_n}(X_i)\,\1_{W_n}(X_j)\,\1_{W_n}(X_k)}{|(W_n-X_i) \cap (W_n - X_j)|}\,\1_{c_n^{\a}\,r\,B}(X_j-X_i)\n\\
&&\phantom{00}\n\\
&+& \ew \Su_{i,k \ge 1} \frac{\1_{W_n}(X_i)\,\1_{W_n}(X_k)\,\1_{c_n^{\a}\,r\,B}(X_k-X_i)}{|(W_n-X_i) \cap (W_n - X_k)|}
+ \ew \Su_{j,k \ge 1} \frac{\1_{W_n}(X_j)\,\1_{W_n}(X_k)\,\1_{c_n^{\a}\,r\,B}(X_k-X_j)}{|(W_n-X_k) \cap (W_n - X_j)|}\n\\
&&\phantom{00}\n\\
&=& \int_{\R^d}\int_{\R^d}\int_{\R^d}\frac{\1_{W_n}(x)\,\1_{W_n}(y)\,\1_{W_n}(z)}{|W_n \cap (W_n - y + x)|}\,\1_{c_n^{\a}\,r\,B}(y-x)\,
\a^{(3)}({\rm d}(x,y,z))
\n\eea

\bea 
&+& 2\,\int_{\R^d}\int_{\R^d}\frac{\1_{W_n}(x)\,\1_{W_n}(y)}{|W_n \cap (W_n - y + x)|}\,\1_{c_n^{\a}\,r\,B}(y-x)\,\a^{(2)}({\rm d}(x,y))\n\\
&&\phantom{00}\n\\
&=& \l\,\int_{\R^d}\int_{\R^d}\int_{\R^d}\frac{\1_{W_n}(x)\,\1_{W_n}(y+x)\,\1_{W_n}(z+x)}{|W_n \cap (W_n - y)|}\,\1_{c_n^{\a}\,r\,B}(y)\,
\a_{red}^{(3)}({\rm d}(y,z))\,{\rm d}x\n\\
&&\phantom{00}\n\\
&+& 2\,\l\,\int_{\R^d}\int_{\R^d}\frac{\1_{W_n}(x)\,\1_{W_n}(y+x)}{|W_n \cap (W_n - y)|}\,\1_{c_n^{\a}\,r\,B}(y)\,
\a_{red}^{(2)}({\rm d}y)\,{\rm d}x\n\\
&&\phantom{00}\n\\
&=& \l\,\int_{\R^d}\int_{\R^d}\frac{|W_n \cap (W_n-y) \cap (W_n-z)|\,\,\1_{c_n^{\a}\,r\,B}(y)}{|W_n \cap (W_n - y)|}\,
\a_{red}^{(3)}({\rm d}(y,z)) + 2\,\l\,\a_{red}^{(2)}(c^{\a}_n\,r\,B)\,.
\label{T4}\eea

\medskip
Disintegration of formula (\ref{ag}) for $k=3$ on both sides  w.r.t. the third component leads to 

\bea
\a_{red}^{(3)}(A_1 \times A_2)&=&\g_{red}^{(3)}(A_1 \times A_2) +  \l\,|A_1|\,\g_{red}^{(2)}(A_2)
+ \g^{(2)}(A_1\times A_2) + \l\,\g_{red}^{(2)}(A_1)\,|A_2| + \l^2\,|A_1|\,|A_2|\n\\
&&\phantom{00}\n\\
&=&\g_{red}^{(3)}(A_1 \times A_2)  + \l\,|A_1|\,\g_{red}^{(2)}(A_2) + \g^{(2)}(A_1\times A_2) + \l\,\a_{red}^{(2)}(A_1)\,|A_2| 
\label{3red}\eea

\mn
for any bounded $A_1, A_2 \in {\mathcal B}^d\,$. Thus, according to (\ref{3red}) we rewrite (\ref{T4}) so that 
$\ew N(W_n)\,(\widehat{\l^2 K_B})_n(c_n^{\a}\,r)$  takes the form 

\bea
&&\l\,\int_{\R^d}\int_{\R^d}\frac{|W_n \cap (W_n-y) \cap (W_n-z)|\,\1_{c_n^{\a} r B}(y)}{|W_n \cap (W_n - y)|}\,
\big[ \g_{red}^{(3)}({\rm d}(y,z)) + \l\,\g_{red}^{(2)}({\rm d}z)\,{\rm d}y + \g^{(2)}({\rm d}(y,z))\big]\n\\
&&\phantom{00}\n\\
&+&  \l^2\,\int_{\R^d}\int_{\R^d}\frac{|W_n \cap (W_n-y) \cap (W_n-z)|\,\1_{c_n^{\a} r B}(y)}{|W_n \cap (W_n - y)|}\,
\a_{red}^{(2)}({\rm d}y)\,{\rm d}z + 2\,\l\,\a_{red}^{(2)}(c^{\a}_n\,r\,B)\,.
\n\eea                               

Since $\int_{\R^d}|W_n \cap (W_n-y) \cap (W_n-z)|\,{\rm d}z = |W_n \cap (W_n-y)|\,|W_n|$
 multiple application of Fubini's theorem  yields the estimate 

\be 
\Big|\,\ew N(W_n)\,(\widehat{\l^2 K_B})_n(c_n^{\a}\,r) - \l^2\,|W_n|\,\a_{red}^{(2)}(c_n^{\a}\,r\,B) 
- 2\,\l\,\a_{red}^{(2)}(c^{\a}_n\,r\,B) - 2\,\l^2\,|c_n^{\a}\,r\,B|\,\g_{red}^{(2)}(\R^d)\,\Big| 
\n\ee
\be \le  \l\,\|\g_{red}^{(3)}\|_{var} + \l\,\int_{\R^d}\int_{c_n^{\a}\,r\,B} \Bigl(1 - \frac{|W_n \cap (W_n-y) 
\cap (W_n-z)|}{|W_n \cap (W_n - y)|}\Bigr)\,{\mathrm d}y \,\g_{red}^{(2)}({\mathrm d}z)\,.
\n\ee

\medskip
From Lemma 1 we see that, for any $y \in c_n^{\a}\,r\,B$ and $z \in \R^d$,

\be 
1 - \frac{|W_n \cap (W_n-y) \cap (W_n-z)|}{|W_n|} \le 1 - \frac{|W_n \cap (W_n - y)|}{|W_n|} + 
1 - \frac{|W_n \cap (W_n - z)|}{|W_n|} \longn 0\,,
\n\ee

\noindent
whence together with (\ref{ared}) and Lebesgue's convergence theorem  it follows (\ref{rest}) so that (\ref{T3}) is proved. 

In order to verify the second limit of (\ref{T12}) we use again the relations $\ew (\widehat{\l^2 K_B})_n(c_n^{\a}\,r) 
= \l\,\a_{red}^{(2)}(c^{\a}_n\,r\,B)$ and $\ew N(W_n) = \l\,|W_n|$ which allows to write $T^{(2)}_n(r\,B)$ in the following form:

\bea
T^{(2)}_n(r\,B) &=& \ew \Bigl[\,|W_n|\,(\widehat{\l^2 K_B})_n(c_n^{\a}\,r) - N(W_n)\,\a_{red}^{(2)}(c^{\a}_n\,r\,B)\,\Bigr]^2\n\\
&&\phantom{00000000}\n\\
&=& |W_n|^2\,\ew \bigl[(\widehat{\l^2 K_B})_n(c_n^{\a}\,r)\bigr]^2 - \bigl(\a_{red}^{(2)}(c^{\a}_n\,r\,B)\bigr)^2\,\ew N^2(W_n) 
- T^{(3)}_n(r\,B)\n\\
&&\phantom{00000000}\n\\
&=& |W_n|^2\,\var \bigl[(\widehat{\l^2 K_B})_n(c_n^{\a}\,r)\bigr] - \var\big(N(W_n))\,\bigl(\a_{red}^{(2)}(c^{\a}_n\,r\,B)\bigr)^2 
- T^{(3)}_n(r\,B)\,.
\n\eea

Having in mind the just proved limit (\ref{T3}) and the limits  ${\var}\big(N(W_n))/|W_n| \longn 
\l\,\big(1+\g_{red}^{(2)}(\R^d)\big)$ and (\ref{ared}), we need to show that
\be 
\lim_{n \to \infty} |W_n|^{1-2\,\a}\,{\var}\bigl[(\widehat{\l^2 K_B})_n(c_n^{\a}\,r)\bigr]
 = 4\,\l^3\,|r\,B|^2\,\bigl(1+\g_{red}^{(2)}(\R^d)\bigr)\,. 
\label{varKdach}\ee

According to the definition of the estimator $(\widehat{\l^2 K_B})_n(r)$ we first rewrite its second moment as expectations of 
multiple sums over pairwise distinct atoms of $N$ which can be expressed  in terms
of integrals w.r.t. the factorial moment measures $\a^{(k)}(\cdot)$ for $k=2,3,4\,$. In this way we get

\bea
\ew \bigl[(\widehat{\l^2 K_B})_n(r)\bigr]^2 = \;{\ew}\sum_{i,j,k,\ell \ge 1 \atop i \ne j, k \ne \ell}\frac{\1_{W_n}(X_i)\,\1_{W_n}(X_j)\,
\1_{W_n}(X_k)\,\1_{W_n}(X_\ell)\,\1_{r\,B}(X_j-X_i)\,\1_{r\,B}(X_\ell-X_k)}{|(W_n-X_i) \cap (W_n - X_j)|\,|(W_n-X_k) \cap (W_n - X_\ell)|}&&\n\\
\phantom{00000000}&&\n\\
= \;\ew \Su_{i,j,k,\ell \ge 1}\frac{\1_{W_n}(X_i)\,\1_{W_n}(X_j)\,
\1_{W_n}(X_k)\,\1_{W_n}(X_\ell)\,\1_{r\,B}(X_j-X_i)\,\1_{r\,B}(X_\ell-X_k)}{|(W_n-X_i) \cap (W_n - X_j)|\,|(W_n-X_k) \cap (W_n - X_\ell)|}&&\n\\
\phantom{00000000}&&\n\\
+ \;4\;\ew \Su_{i,j,k \ge 1}\frac{\1_{W_n}(X_i)\,\1_{W_n}(X_j)\,\1_{W_n}(X_k)\,\1_{r\,B}(X_j-X_i)\,\1_{r\,B}(X_k-X_i)}{|(W_n-X_i) \cap (W_n - X_j)|\,
|(W_n-X_i) \cap (W_n - X_k)|}
\phantom{00000000}&&\n\\
\phantom{00000000}&&\n\\
+ \;2\;\ew \Su_{i,j\ge 1}\frac{\1_{W_n}(X_i)\,\1_{W_n}(X_j)\,\1_{r\,B}(X_j-X_i)}{|(W_n-X_i) \cap (W_n - X_j)|^2}\,,\phantom{00000000000000000000}&&
\n\eea

\be 
\bigl[\ew (\widehat{\l^2 K_B})_n(r)\bigr]^2 = \int_{W_n}\int_{W_n}\int_{W_n}\int_{W_n} \frac{\1_{r\,B}(y - x)\,\1_{r\,B}(v - u)\,
\a^{(2)}({\mathrm d}(x,y))\,\a^{(2)}({\mathrm d}(u,v))}{|W_n \cap (W_n-y+x)|\,|W_n \cap (W_n-v+u)|}
\n\ee

and
\bea
\var \bigl[(\widehat{\l^2 K_B})_n(c_n^{\a}\,r)\bigr] &=& \int_{\R^d}\int_{\R^d}\int_{\R^d}\int_{\R^d} \frac{\1_{W_n}(x)
\1_{W_n}(y)\1_{W_n}(u)\1_{W_n}(v)\,\1_{c_n^{\a}\,r\,B}(y - x)}{|W_n \cap (W_n-y+x)|\,|W_n \cap (W_n-v+u)|}\n\\
&\times& \1_{c_n^{\a}\,r\,B}(v - u)\,\big[\,\a^{(4)}({\mathrm d}(x,y,u,v)) - \a^{(2)}({\mathrm d}(x,y))\,\a^{(2)}({\mathrm d}(u,v))\,\big]  \n\\
&&\phantom{0000}\n\\
&+& \;4\;\int_{W_n}\int_{W_n}\int_{W_n}\frac{\1_{c_n^{\a}\,r\,B}(y - x) \,\1_{c_n^{\a}\,r\,B}(z - x)\,
\a^{(3)}({\mathrm d}(x,y,z))}{|W_n \cap (W_n-y+x)|\,|W_n \cap (W_n-z+x)|}\n\\
&&\phantom{0000}\n\\
&+& \;2\;\int_{W_n}\int_{W_n} \frac{\1_{c_n^{\a}\,r\,B}(y - x)\,\a^{(2)}({\mathrm d}(x,y)) }{|W_n \cap (W_n-y+x)|^2} \n\\
&&\phantom{0000}\n\\
&=& T^{(4)}_n(r\,B) + 4\,T^{(5)}_n(r\,B) + 2\,T^{(6)}_n(r\,B)\;.
\label{VarKBcnr}\eea

\medskip
We first treat the asymptotic behaviour of $T^{(5)}_n(r\,B)$ and $T^{(6)}_n(r\,B)$. After some obvious rearrangements combined 
with Lemma 1 and (\ref{ared}) we obtain  

\bea
\frac{|W_n|\,T^{(6)}_n(r\,B)}{|W_n|^{2\,\a}} &=& \frac{\l\,|W_n|}{|W_n|^{2\,\a}}\,\int_{\R^d}\int_{\R^d} \frac{\1_{W_n}(x)\,\1_{W_n}(y+x)\,
\1_{c_n^{\a}\,r\,B}(y)}{|W_n \cap (W_n-y)|^2}\,\a^{(2)}_{red}({\mathrm d}y){\mathrm d}x\n\\
&&\phantom{0000}\n\\
&=& \frac{\l\,|W_n|}{|W_n|^{2\,\a}}\,\int_{\R^d} \frac{\1_{c_n^{\a}\,r\,B}(y)\,\a^{(2)}_{red}({\mathrm d}y)}{|W_n \cap (W_n-y)|}\quad 
(\;\longn \;\; \l\,\a^{(2)}_{red}(r\,B)\;\;\mbox{for}\;\; \a = 0\;)\n\\
&&\phantom{0000}\n\\
&\le & \l\,\sup_{y\in c_n^{\a}\,r\,B} \frac{|W_n|}{|W_n \cap (W_n-y)|}
\,\frac{\a^{(2)}_{red}(c_n^{\a}\,r\,B)}{|W_n|^{2\,\a}} \longn 0\quad\mbox{for}\;\; \a >  0\,.
\label{T6}\eea

\medskip
Similarly, using the differential reduction formula $\a^{(3)}({\mathrm d}(x,y,z)) = \l\,\a_{red}^{(3)}({\rm d}y -x, {\rm d}z -x)\,{\rm d}x$  and Lemma 1 we find that

\bea 
\frac{|W_n|\,T^{(5)}_n(r\,B)}{|W_n|^{2\,\a}}&=&\l\,\int_{\R^d}\int_{\R^d}\frac{\1_{c_n^{\a}\,r\,B}(y) \,\1_{c_n^{\a}\,r\,B}(z)\,|W_n|\,
|W_n \cap (W_n-y)\cap (W_n-z)|}{|W_n|^{2\,\a}\,|W_n \cap (W_n-y)|\,|W_n \cap (W_n-z)|} \,\a_{red}^{(3)}({\mathrm d}(y,z))
\n\\
&&\phantom{0000}\n\\
&=& \frac{\l\,\a_{red}^{(3)}(c_n^{\a}\,r\,B \times c_n^{\a}\,r\,B)}{|W_n|^{2\,\a}}\,\big( 1 + \t\,\e_n(r)\big)\,,
\label{T5}\eea
where $\t \in [-1,1]$ is suitably chosen and   

\be 
\e_n(r) := \sup_{y,z\in c_n^{\a}\,r\,B}\Bigg|\,\frac{|W_n|\,|W_n \cap (W_n-y)\cap (W_n-z)|}{|W_n \cap (W_n-y)|\,|W_n \cap (W_n-z)|} - 1\,\Bigg| \longn 0
\quad\mbox{for any fixed}\;\;r > 0\,,
\n\ee

\mn
which shows that $|W_n|\,T^{(5)}_n(r\,B) \longn \l\,\a_{red}^{(3)}(r\,B \times r\,B)$ for $\a = 0\,$. 
On the other hand, by (\ref{3red}), 
\bea 
\l\,\a_{red}^{(3)}(c_n^{\a}\,r\,B \times c_n^{\a}\,r\,B) &=& \l^2\,\a_{red}^{(2)}(c_n^{\a}\,r\,B)\,| c_n^{\a}\,r\,B| + \l^2\,\g_{red}^{(2)}(c_n^{\a}\,r\,B)\,
| c_n^{\a}\,r\,B|\n\\
&+& \l\,\g_{red}^{(3)}(c_n^{\a}\,r\,B \times c_n^{\a}\,r\,B) + \l^2\,\int_{\R^d} | c_n^{\a}\,r\,B \cap (c_n^{\a}\,r\,B-y) |\,\g_{red}^{(2)}({\mathrm d}y)\,.   
\n\eea

\mn
Hence, together with (\ref{ared}) it is easily seen that $\quad |W_n|^{1-2\,\a}\,T^{(5)}_n(r\,B) \longn \l^3\,|r\,B|^2\,$.

\medskip
To prove (\ref{varKdach}) it remains to  verify the limit
\be  
\lim_{n \to\infty} \frac{|W_n|\,T^{(4)}_n(r\,B)}{|W_n|^{2\,\a}} = 4\,\l^3\,|r\,B|^2\, \g_{red}^{(2)}(\R^d)\quad\mbox{for}\;\; \a >  0\,\,.
\label{4Int}\ee

\medskip
For proving this we make use of the decomposition

\bea 
&&\a^{(4)}({\mathrm d}(x,y,u,v))\, - \,\a^{(2)}({\mathrm d}(x,y))\,
\a^{(2)}({\mathrm d}(u,v)) = \g^{(4)}({\mathrm d}(x,y,u,v)) \label{a4a2a2}\\
&&\phantom{00}\n\\
&+& \l\,\big[\,{\rm d}x\,\g^{(3)}({\mathrm d}(y,u,v))\,+\,{\mathrm d}y\,\g^{(3)}({\mathrm d}(x,u,v))\,
+\,{\mathrm d}u\,\g^{(3)}({\mathrm d}(x,y,v))\,+\,{\mathrm d}v\,\g^{(3)}({\mathrm d}(x,y,u))\,\big]\n\\
&&\phantom{00}\n\\
&+& \g^{(2)}({\mathrm d}(x,u))\,\g^{(2)}({\mathrm d}(y,v))\,+\,\g^{(2)}({\mathrm d}(x,v))\,\g^{(2)}({\mathrm d}(y,u))\n\\ 
&&\phantom{00}\n\\ 
&+& \l^2\,\big[\,{\mathrm d}x\,{\mathrm d}u\,\g^{(2)}({\mathrm d}(y,v))+ {\mathrm d}x\,{\mathrm d}v\,\g^{(2)}({\mathrm d}(y,u))\, 
+ {\mathrm d}y\,{\mathrm d}u\,\g^{(2)}({\mathrm d}(x,v))\,
+ {\mathrm d}y\,{\mathrm d}v\,\g^{(2)}({\mathrm d}(x,u))\,\big]\,,
\n\eea

\sn
which is obtained by applying the formula (\ref{ag}) for $k=2$ and $k=4$ with $A_1={\rm d}x$, $A_2={\rm d}y$, 
$A_3={\rm d}u$ and $A_4={\rm d}v\,$. After rewriting the integrals in terms of  reduced cumulant measures 
and in view of the assumptions $\|\g_{red}^{(k)}\|_{var} < \infty$ for $k=2,3,4$ it turns out that, for $\a > 0\,$, 
only the below four integrals contribute with a non-zero limit to the r.h.s. of $(\ref{4Int})$:

\bea
&& \frac{\l^2\,|W_n|}{|W_n|^{2\a}}\,\int_{\R^d}\int_{\R^d}\int_{\R^d}\int_{\R^d} \frac{\1_{W_n}(x)\1_{W_n}(y)\1_{W_n}(u)\1_{W_n}(v)\,
\1_{c_n^{\a}\,r\,B}(y - x)\,\1_{c_n^{\a}\,r\,B}(v - u)}{|W_n \cap (W_n-y+x)|\,|W_n \cap (W_n-v+u)|}\n\\
\phantom{000}\n\\
&& \times \big[\,\,\g^{(2)}({\mathrm d}(y,v))\,{\mathrm d}x\,{\mathrm d}u+ {\mathrm d}x\,{\mathrm d}v\,\g^{(2)}({\mathrm d}(y,u))\,
+ {\mathrm d}y\,{\mathrm d}u\,\g^{(2)}({\mathrm d}(x,v))\,+ {\mathrm d}y\,{\mathrm d}v\,\g^{(2)}({\mathrm d}(x,u))\,\big]\n\\
\phantom{000}\n\\
&=& \frac{4\,\l^3\,|W_n|}{|W_n|^{2\a}}\,\int_{\R^d}\int_{\R^d}\int_{\R^d}\int_{\R^d}\frac{\1_{W_n}(x+y)\,\1_{W_n}(u+v+y)\,
\1_{W_n}(y)\,\1_{W_n}(v+y)}{|W_n \cap (W_n+x)|\,|W_n \cap (W_n+u)|}\n\\
\phantom{000}\n\\
&&\phantom{0000000000000000000000000000}\times \,\1_{c_n^{\a}\,r\,B}(x)\,\1_{c_n^{\a}\,r\,B}(u)\,\g^{(2)}_{red}({\rm d}v)\,{\rm d}x\,{\rm d}u\n\\
\phantom{000}\n\\
&=& 4\,\l^3\,\int_{\R^d}\int_{r\,B}\int_{r\,B}\frac{|W_n|\,|W_n\cap (W_n- c_n^{\a}\,x)\cap (W_n-v)\cap (W_n- c_n^{\a}\,u -v)|}{|W_n \cap 
(W_n- c_n^{\a}\,x)|\,|W_n \cap (W_n- c_n^{\a}\,u)|}\,{\mathrm d}x\,{\mathrm d}u\,\g^{(2)}_{red}({\mathrm d}v)\n\\
\phantom{000}\n\\
&& \longn 4\,\l^3\,|r\,B|^2\,\g^{(2)}_{red}(\R^d)\quad\mbox{for}\;\; \a \ge 0\,,
\n\eea

\mn
where the latter limit is justified by Lemma 1 and and Lebesgue's convergence theorem.
Thus, $(\ref{4Int})$ is  proved which completes the proof of Lemma 2.  $\Box$

\bn
For $\a = 0$ the remaining seven integrals on the r.h.s. of (\ref{a4a2a2}) possesses non-zero limits in general:

\bea
|W_n|\,\int_{\R^d}\int_{\R^d}\int_{\R^d}\int_{\R^d} \frac{\1_{W_n}(x)\1_{W_n}(y)\1_{W_n}(u)\1_{W_n}(v)\,
\1_{r B}(y - x)\,\1_{r B}(v - u)}{|W_n \cap (W_n-y+x)|\,|W_n \cap (W_n-v+u)|}\,\big[\g^{(4)}({\rm d}(x,y,u,v))&& \n\\
\phantom{000}&&\n\\
+ 2\,\l\,\g^{(3)}({\rm d}(y,u,v))\,{\rm d}x 
+ 2\,\l\,\g^{(3)}({\rm d}(x,y,v))\,{\rm d}u + 2\,\g^{(2)}({\rm d}(x,u))\,\g^{(2)}({\rm d}(y,v))\big] \phantom{00000} &&\n\\
 \longn \;\; \l\,\int_{\R^{3d}}\1_{r B}(v-u)\,\1_{r B}(y)\,\g^{(4)}_{red}({\rm d}(y,u,v)) 
+ 2\,\l^2\,|r B|\,\int_{\R^{2d}}\1_{r B}(v-u)\,\g^{(3)}_{red}({\rm d}(u,v))&&\n\\
+ 2\,\l^2\,|r B|\,\g^{(3)}_{red}(r B \times \R^d) 
+ 2\,\l^2\,\int_{\R^{2d}}|(r B -u) \cap (r B -v) |\,\g^{(2)}_{red}({\rm d}u)\,\g^{(2)}_{red}({\rm d}v)\,. \phantom{0000}
\n\eea

\noindent
All above limits obtained for $\a = 0$ and inserted in (\ref{VarKBcnr}) are summarized in 

\bn
\textbf{Remark 7} \ For any  $B_4$-mixing PP  $N = \sum_{i\ge 1} \d_{X_i}$ on $\R^d$ the asymptotic variance 
\be 
\tau^2_B(r) := \lim_{n\to \infty} |W_n|\,\var\bigl(\,(\widehat{\l^2 K_B})_n(r)\,\bigr) 
\n\ee
exists and takes the form

\bea 
\tau^2_B(r) & = &
\l\,\int_{\R^{3\,d}}\1_{r\,B}(y-x)\,\1_{r\,B}(z)\,\g^{(4)}_{red}({\rm d}(x,y,z)) + 4\,\l^2\,|r B|\,\g^{(3)}_{red}(r B \times \R^d)\n\\
&+& 2\,\l^3\,\int_{\R^{2\,d}}|(r\,B-x) \cap (r\,B-y) |\,\g^{(2)}_{red}({\rm d}x)\,\g^{(2)}_{red}({\rm d}y)
+ 4\,\l^3\,|r\,B|^2\,\g^{(2)}_{red}(\R^d)\n\\
&+& 4\,\l\,\a^{(3)}_{red}(r\,B \times r\,B)
+2\,\l\,\a^{(2)}_{red}(r\,B)\,.
\label{asvar}\eea

\bn
{\em Proof of Lemma 3} \ For notational ease  we put

\be 
(\widehat{\l \a_{red}^{(2)}})_n(A) :=  \Su_{i,j \ge 1}\frac{\1_{W_n}(X_i)\,\1_{W_n}(X_j)\,
\1_{c_n^{\a}\,A}(X_j-X_i)}{|(W_n - X_i) \cap (W_n - X_j)|}\quad\mbox{for bounded}\;\; A \in {\mathcal B}^d\,.
\n\ee

\medskip
From (\ref{Dna}) and the definition of $\a_{red}^{(2)}(\cdot)$ we get, for $0 \le s \le t$,

\be 
\D_{B,n}^{(\a)}(t)-\D_{B,n}^{(\a)}(s) = |W_n|^{1/2-\a}\,(\widehat{\l\,\a_{red}^{(2)}})_n((t B) \setminus (s B))\quad\mbox{and}\quad 
\ew (\widehat{\l \a_{red}^{(2)}})_n(A) = \l\,\a_{red}^{(2)}(c^{\a}_n\,A)\,.
\n\ee

Therefore, it suffices to show that

\be 
|W_n|^{1-2\,\a}\,\var\bigl(\,(\widehat{\l \a_{red}^{(2)}})_n(A)\,\bigr) \le a_1\,|A|^2 + a_2\,|A|\,|W_n|^{-\a} + a_3\,|W_n|^{-2\,\a}
\label{A}\ee

\mn
for any bounded, $\O$-symmetric ( i.e. $A = - A$) $A \in {\mathcal B}^d\,$.   

\medskip 
For this purpose we use the decomposition (\ref{VarKBcnr}) and replace the ball $r\,B$ by a bounded $A \in {\mathcal B}^d$. Due to Lemma 1 
we find some $n_0$ (depending on $\sup\{\|x\| : x \in A\}$) such that $\sup_{y \in c_n^{\a}\,A} |W_n|/|W_n \cap (W_n-y)| \le 2$  for 
all $n \ge n_0$. From (\ref{T6}) with $A$ instead of $r\,B$  together with (\ref{density}) we arrive at

\be 
|W_n|^{1-2\,\a}\,T^{(6)}_n(A) \le 2\,\l\,\a^{(2)}_{red}(c_n^{\a}\,A)\,|W_n|^{-2\,\a} \le 2\,a_0\,\l\,|A|\,|W_n|^{-\a}
\quad\mbox{for}\;\;n \ge n_0\,.
\label{T6A}\ee

\medskip
In the same way we modify (\ref{T5}), where $1 + \t\,\e_n(r)$ can be replaced by 2 for $n \ge n_0$ (due to the argument used before), 
and obtain  
 
\be 
|W_n|^{1-2\,\a}\,T^{(5)}_n(A) \le 2\,\l\,\a_{red}^{(3)}(c_n^{\a}\,A\times c_n^{\a}\,A)\,|W_n|^{-2\,\a}\quad\mbox{for}\;\;n \ge n_0\,.
\n\ee

Condition (\ref{density}) applied to (\ref{3red}) leads to

\be 
\a_{red}^{(3)}(A_1 \times A_2) \le \|\g_{red}^{(3)}\|_{var} +  3\,\l\,a_0\,|A_1|\,|A_2| - 2\,\l^2\,|A_1|\,|A_2|\,,
\n\ee
whence it follows that

\be 
|W_n|^{1-2\,\a}\,T^{(5)}_n(A) \le 2\,\l\,\|\g_{red}^{(3)}\|_{var}\,|W_n|^{-2\,\a} + 6\,\l^2\,a_0\,|A|^2\quad\mbox{for}\;\;n \ge n_0\,.
\label{T5A}\ee

\medskip
On the other hand, (\ref{3red})and the first condition of (\ref{density34}) yield $\a_{red}^{(3)}(A_1 \times A_2) \le \max\{a_{01}+2\,\l^2, 
3\,\l\,a_0\}\,|A_1|\,|A_2|$ so that $|W_n|^{1-2\,\a}\,T^{(5)}_n(A) \le a_1\,|A|^2$ for $n \ge n_0$ for some constant $a_1 > 0$.

\bigskip
To obtain suitable bounds of $T^{(4)}_n(A)$ (defined  in (\ref{VarKBcnr}) for $A = r\,B$ ) for any bounded, $\O$-symmetric set 
$A \in {\mathcal B}^d$ we have to estimate all integrals w.r.t. the measures on the r.h.s of
(\ref{a4a2a2}). As before we choose $n_0$ such that $|W_n| \le 2\,\inf_{x \in c_n^{\a}\,A}|W_n \cap (W_n-x)|$ for $n \ge n_0$.
By the $\O$-symmetry of $A$ it is easily checked that  

\bea
T^{(4)}_n(A) = \int_{\R^d}\int_{\R^d}\int_{\R^d}\int_{\R^d} \frac{\1_{W_n}(x)\1_{W_n}(y)\1_{W_n}(u)\1_{W_n}(v)\,
\1_{c_n^{\a}\,A}(y - x)\,\1_{c_n^{\a}\,A}(v - u)}{|W_n \cap (W_n-y+x)|\,|W_n \cap (W_n-v+u)|}\,\times \phantom{00000}&&\n\\
\phantom{000}&&\n\\
\;\big[\g^{(4)}({\rm d}(x,y,u,v)) + 4\,\l\,\g^{(3)}({\mathrm d}(y,u,v)){\rm d}x + 2\,\g^{(2)}({\rm d}(x,u))\,\g^{(2)}({\rm d}(y,v)) + 
4\,\l^2\,\g^{(2)}({\mathrm d}(y,v)) {\mathrm d}(x,u)\big].
\n\eea

\medskip
After passing to the  reduced factorial cumulant measures $\g_{red}^{(k)}(\cdot)$ and using that  $\|\g_{red}^{(k)}\|_{var} < \infty$ 
for $k=2,3,4$ we obtain after some straightforward, but lengthy  calculations the estimate

\be 
|W_n|\,\frac{T^{(4)}_n(A)}{|W_n|^{2\,\a}} \,\le\, 2\,\l\, \frac{\|\g_{red}^{(4)}\|_{var}}{|W_n|^{2\,\a}}+ 
8\,\l^2\,|A|\,\frac{\|\g_{red}^{(3)}\|_{var}}{|W_n|^{\a}} + 8\,\l^2\,|A|\,\frac{\|\g_{red}^{(2)}\|_{var}^2}{|W_n|^{\a}} + 
8\,\l^3\,|A|^2\,\|\g_{red}^{(2)}\|_{var}
\label{T4A}\ee
for all $n \ge n_0\,$.  Clearly, (\ref{density}) implies that the reduced covariance measure $\g_{red}^{(2)}(\cdot) = \a_{red}^{(2)}(\cdot)-\l\,|\cdot|$ 
admits the estimate $|\,\g_{red}^{(2)}(A)\,|  \le \max\{a_0,\l\}\,|A|$ for all bounded $A \in B^d$. After  inserting the latter on the r.h.s. 
of (\ref{T4A}) we obtain together with (\ref{T6A}), (\ref{T5A}) the estimate (\ref{A}). Setting $A = (t B) \setminus (s B)$ in (\ref{A} terminates 
the proof of Lemma 3. $\Box$

\bigskip
To prove (\ref{rem1}) we use $\g^{(4)}({\mathrm d}(x,y,u,v)) = \l\,\kappa^{(4)}(y-x,u-x,v-x)\,{\mathrm d}(v,u,y,x)$ and (\ref{density34}) 
to find an appropriate bound of 

\bea 
&& \frac{|W_n|}{|W_n|^{2\a}}\,\Bigg| \int_{\R^{4\,d}}\frac{\1_{W_n}(x)\1_{W_n}(y)\1_{W_n}(u)\1_{W_n}(v)\,
\1_{c_n^{\a}\,A}(y - x)\,\1_{c_n^{\a}\,A}(v - u)}{|W_n \cap (W_n-y+x)|\,|W_n \cap (W_n-v+u)|}\,\g^{(4)}({\mathrm d}(x,y,u,v))\Bigg|\n\\
&& \phantom{00}\n\\
&\le & \frac{|W_n|}{|W_n|^{2\a}}\, \int_{\R^{3\,d}}\frac{|W_n \cap (W_n-y) \cap (W_n-u) \cap (W_n-v+u)|\,
\1_{c_n^{\a}\,A}(y)\,\1_{c_n^{\a}\,A}(v)}{|W_n \cap (W_n-y)|\,|W_n \cap (W_n-v)|}\,\n\\
&& \phantom{00000000000000000000000000000000000000000000000000} \times |\kappa^{(4)}(y,u,v+u)|\,{\mathrm d}(v,u,y)\n\\
&\le &\frac{2}{|W_n|^{2\a}}\,\int_{\R^{2\,d}}\1_{c_n^{\a}\,A}(y)\,\1_{c_n^{\a}\,A}(v)\,\int_{\R^d}|\kappa^{(4)}(y,u,v+u)|\,{\rm d}u\,{\rm d}(v,y)\,.
\n\eea

\medskip
An application of the Cauchy-Schwarz inequality combined with the second condition of (\ref{density34})  yields the following bound of the latter line: 

\be 
\frac{2\,|c_n^{\a}\,A|}{|W_n|^{2\a}}\,\left(\int_{\R^{2\,d}} \Big(\int_{\R^d}|\kappa^{(4)}(y,u,v+u)|\,{\rm d}u \Big)^2\,{\rm d}(v,y)\right)^{1/2} 
\le 2\,\sqrt{a_{02}}\,\frac{|A|}{|W_n|^{\a}}\,.
\n\ee
Together with the foregoing estimates of the other terms of $T^{(4)}_n(A)$ we obtain estimate (\ref{A}) without  $a_3\,|W_n|^{-2\,\a}$ but changed 
constants $a_1$ and $a_2$. For $A = (t B) \setminus (s B)$ this coincides with the desired relation  (\ref{rem1}). $\Box$

\bn
\section{Appendix A}

\bigskip 
In this appendix we state and prove a result from convex geometry which has 
relevance for the statistical analysis of spatial PPes and random closed sets 
observed on increasing compact convex sampling  windows. Unfortunately, 
we could not find an appropriate reference to the below estimates
 (\ref{inequ}) in the huge literature on convex geometry.   

\mn  
Let $K \ne \emptyset$  be a {\em convex body} (compact convex set) in $\R^d$ 
with inball radius $\vr(K) > 0$. In this case the corresponding {\em 
quermassintegrals} $W_k(K)$ are also positive for $k = 0,1,...,d\,$, where
$W_0(\cdot) = |\cdot|$ and $d\,W_1(\cdot) = {\mathcal H}^{d-1}(\pa(\cdot))\,$, 
stands for the $d$-dimensional volume and the $(d-1)$-dimensional surface content,
respectively. Let denote $K \oplus (r B_e) := \{ x \in \R^d : 
(r B_e) + x \cap K \ne \emptyset \}$ resp.  $K \ominus (r B_e):= \{x \in K:
(r B_e) + x \subseteq K\}$ the dilation resp. erosion of $K$ by a closed 
Euclidean ball $r B_e$ with radius $r \ge 0$ centered at the origin $\O$. The following 
lemma gives bounds of the deviation of the volumes $|K \oplus (r B_e)|$ 
and $|K \ominus (r B_e)|$ from $|K|$ in terms of the ratio $r/\vr(K)\,$. 

\bn
\textbf{Lemma 5} \ With the above notation the following two inequalities 
hold for $0 \le r \le \vr(K)\,$:
\be
\label{inequ}\qquad
0 \le 1 - \frac{|K \ominus (r B_e)|}{| K |} \le \frac{d\,r}{\vr(K)} 
\qquad\mbox{{\em and}}\qquad\frac{r}{\varrho(K)} \le 
\frac{|K \oplus (r B_e)|}{| K |} - 1 \le \frac{(2^d-1)\,r}{\vr(K)}\,.
\ee

\bn
{\em Proof of Lemma 5} \ We start by recalling the fact that the function $r \mapsto |K \ominus r B_e|$ 
has a continuous derivative ${\mathcal H}^{d-1}(\pa(K \ominus (r B_e))$ for $0 \le r < \vr(K)\,$, 
see Hadwiger (1957), p. 207. This fact and the isotony of the quermassintegral $W_1(\cdot)$ yield 

\be 
0 \le |K| - |K \ominus (r B_e)| = \int_0^r {\mathcal H}^{d-1}(\pa(K\ominus(\rho B_e))){\mathrm d}\rho 
\le r\,{\mathcal H}^{d-1}(\pa K)\;\;\mbox{for}\;\; 0 \le r \le \vr(K)\,.
\label{omi}\ee

\mn
Since $|K \ominus (\vr(K) B_e)| = 0$ we get immediately that $|K| \le \vr(K)\,{\mathcal H}^{d-1}(\pa K)\,$. 
On the other hand, it was proved in Wills (1970) that $d\,|K| \ge \vr(K)\,{\mathcal H}^{d-1}(\pa K)$ which  
implies the inclusion

\be 
\frac{1}{\vr(K)} \le \frac{{\mathcal H}^{d-1}(\pa K)}{|K|} \le \frac{d}{\vr(K)}\;.
\label{rat}\ee

\mn
Combining (\ref{omi}) and (\ref{rat}) yields the first inequality of the
Lemma. To verify the second assertion we make use of the famous
Steiner formula, see Hadwiger (1957), which reads as follows:

\be 
\label{ste}
r\,{\mathcal H}^{d-1}(\partial K) = r\,d\,W_1(K) \le |K \oplus (r B_e)| - |K| 
= \sum_{k=1}^d\,{d \choose k}\, W_k(K)\,r^k
\ee

\noindent
Now, we apply of the well-known inequality  $W_{k-1}(K)\,W_{k+1}(K) \le (W_k(K))^2$ for $k = 1,\ldots,d-1\,$, 
see Hadwiger (1957), p. 282, which implies together with the right-hand inequality of (\ref{rat}) that  

\be
W_k(K) \le W_0(K)\,\left(\,\frac{W_1(K)}{W_0(K)}\,\right)^k = 
|K|\,\left(\,\frac{{\mathcal H}^{d-1}(\pa K)}{d\,|K|}\,\right)^k  \le
\frac{|K|}{(\vr(K))^k}\quad\mbox{for}\quad k=1,\ldots,d\;.   
\n\ee

\mn
Inserting this estimate on the r.h.s. of (\ref{ste}) we arrive at
\be
|K \oplus (r B_e)| - |K| \le |K|\,\bigg(\Big(\,1 + \frac{r}{\vr(K)}\,\Big)^d - 1\,\bigg) 
= \frac{r\,|K|}{\vr(K)}\,\sum_{k=0}^{d-1}\,\Bigl(\,1 + \frac{r}{\vr(K)}\,\Bigr)^k \,.
\n\ee
Since $0 \le r/\varrho(K) \le 1\,$, the latter sum is bounded by $2^d - 1$. 
This and the previous estimate combined with the lower bounds in
(\ref{rat}) and (\ref{ste}) yield the second estimate of 
(\ref{inequ}) which completes the proof the Lemma 5. $\Box$

\bn
\textbf{Corollary 4} \ It holds that $|\partial K \oplus (r B_e)| \le (2^d-1+d)\,r\,{\mathcal H}^{d-1}(\partial K)$ for $r \ge 0\,$. 

\mn
{\em Proof of Corollary 4} \ Combining both inequalities of Lemma 5 leads to 
\be
|\partial K \oplus (r B_e)| = |K \oplus (r B_e)|-|K| + |K|-|K \ominus (r B_e)|\; 
\le \; (2^d-1+d)\;r\;\frac{|K|}{\vr(K)}  
\n\ee
and the l.h.s. of (\ref{rat}) confirms the assertion of Corollary 4. $\Box$

\vfill\eject

\noindent
\section{Appendix B} 
In this second appendix we study the asymptotic behaviour of two further estimators of $\l^2\,K_B(r)$ which are 
slightly different from (\ref{KBest}). We consider their scaled and as unscaled version separately. To start with 
we define a (so-called naive) estimator  $(\widetilde {\l^2 K_B})_n(r)$ which also turns out unbiased and 
is easier to calculate, 

\be 
(\widetilde{\l^2 K_B})_n(r) := \frac{1}{|W_n|} \sum_{i\ge 1} \1_{W_n}(X_i)\,(N - \d_{X_i})(r\,B + X_i ) =
\frac{1}{|W_n|} \Su_{i,j \ge 1} \1_{W_n}(X_i)\,\1_{r\,B}( X_j - X_i )\,.
\n\ee
However, in order to calculate $(\widetilde {\l^2 K_B})_n(c_n^{\a}\,r)$ for $0 \le r \le R$ and $0\le \a < 1$ the PP 
$N = \sum_{i\ge 1} \delta_{X_i}$ must be observable in the larger window $W_n \oplus (c_n^{\a}\,R\,B)$ or the original 
window $W_n$ must be replaced by the eroded window $W_n \ominus (c_n^{\a}\,R\,B)$ (which means  `minus-sampling'). To avoid 
minus-sampling and as an further alternative to $(\widetilde{\l^2 K_B})_n(r)$ or $(\widehat{\l^2 K_B})_n(r)$ the estimator

\be 
(\overline{\l^2 K_B})_n(r) := \frac{1}{|W_n|} \Su_{i,j \ge 1} \1_{W_n}(X_i)\,\1_{W_n}(X_j)\,\1_{r\,B}( X_j - X_i )
\n\ee
is sometimes used although it is no longer unbiased but at least still asymptotically unbiased in the following sense: 
Provided that $c_n^{\a}/\vr(W_n) \longn 0$ for some $0 \le \a < 1$ it follows by means of Lemma 1 that, as $n \to \infty\,$, 

\be
\ew (\overline{\l^2 K_B})_n(c_n^{\a}\,r) =  \l\,\int_{c_n^{\a}\,r B}\frac{|W_n \cap (W_n - x)|}{|W_n|}\,
\a^{(2)}_{red}({\rm d}x)                 
= \l^2\,K_B(c_n^{\a}\,r)\,\Big( 1 + {\mathcal O}\Big(\frac{c_n^{\a}}{\vr(W_n)}\Big) \Big)\,.
\n\ee

\bn
\textbf{Lemma 1B} \ Let  $N = \sum_{i\ge 1} \d_{X_i}$ be a $B_4$-mixing PP on $\R^d\,$. Further, let $(W_n)$ be a CAS in $\R^d$ 
with inball radius $\vr(W_n)$ and $0 \le \a < 1\,$. If $c_n^{\a}/\vr(W_n) \longn 0\,$, then

\be
\lim_{n \to \infty}\;|W_n|^{1-2\a}\,\ew\Bigl[\,(\widehat{\l^2 K_B})_n(c_n^{\a}\,r) - (\widetilde{\l^2 K_B})_n(c_n^{\a}\,r)\,\Bigr]^2\,=\,0
\n\ee 
and 

\be
\lim_{n \to \infty}\;|W_n|^{1-2\a}\;\var\Bigl[\,(\overline{\l^2 K_B})_n(c_n^{\a}\,r) - (\widetilde {\l^2 K_B})_n(c_n^{\a}\,r)\,\Bigr]^2 \, = \,0\,.
\n\ee 

\mn
{\em Proof of Lemma 1B} \ We first express the variances of $(\widetilde{\l^2 K_B})_n(c_n^{\a}\,r)$ and  $(\overline{\l^2 K_B})_n(c_n^{\a}\,r)$ 
in terms of factorial moment and cumulant measures  of order $k=2, 3, 4\,$. For this we apply general formula for the covariance 
$C(g,h) := \cov\big(\Su_{i,j\ge 1} g(X_i,X_j),\Su_{k,\ell\ge 1} h(X_k,X_\ell)\big)$ which has been proved in Heinrich (1988), p. 97. 
To avoid ambiguities we present this formula once more for fourth-order stationary PP's and Borel-measurable functions $g, h$ on $\R^d\times \R^d$ 
such that the integrals on the r.h.s. exist:

\bea
C(g,h) &=& \int_{\R^{4d}} g(x,y)\,h(u,v)\,\big[\,\a^{(4)}({\rm d}(x,y,u,v)) - \a^{(2)}({\rm d}(x,y))\,
\a^{(2)}({\rm d}(u,v))\big]\label{cgh}\\
&+& \int_{\R^{3d}} g(x,y)\,\big[\,h(x,u)+h(y,u)+h(u,x)+h(u,y)\,\big]\,\a^{(3)}({\rm d}(x,y,u)) \n\\
&+& \int_{\R^{2d}} g(x,y)\,\big[\,h(x,y)+h(y,x)\,\big]\,\a^{(2)}({\rm d}(x,y))\,,
\n\eea
where the signed measure $\a^{(4)}-\a^{(2)}\times\a^{(2)}$ in (\ref{cgh}) allows the decomposition (\ref{a4a2a2}).

We use (\ref{cgh}) together with (\ref{a4a2a2})
for $g(x,y)= h(x,y)= |W_n|^{-1}\,\1_{W_n}(x)\,\1_{c_n^{\a}\,r\,B}(y - x)$ and get

\bea
|W_n|^{1-2\a}\;\var\big[(\widetilde{\l^2 K_B})_n(c_n^{\a}\,r)\big]
\phantom{0000000000000000000000000000000000000000000000000000000000000000000000000}&& \n\\
&&\n\\
= \frac{1}{|W_n|^{1+2\a}}\;\bigg[\int\limits_{\R^{4d}}\1_{W_n}(x)\,\1_{W_n}(u)\,\1_{c_n^{\a}\,r\,B}(y - x)\,\1_{c_n^{\a}\,r\,B}(v - u)\,
\big(\,\a^{(4)} - \a^{(2)} \times \a^{(2)}\,\big)({\rm d}(x,y,u,v))\phantom{00000000000000000000000000000000000000000000} &&\n\\
\phantom{0000000000}&&\n\\
+ \int\limits_{\R^{3d}}\Big(\1_{W_n}(x)\big(1+\1_{W_n}(u) + \1_{W_n}(y)\big) + \1_{W_n}(u)\1_{W_n}(y)\Big)
 \1_{c_n^{\a} r B}(y - x)\,\1_{c_n^{\a} r B}(u - x)\a^{(3)}({\rm d}(x,y,u))\phantom{000000000000000000000000000000000000000000}&& \n\\
\phantom{000000000}&&\n\\
+ \int\limits_{\R^{2d}} \1_{W_n}(x)\,\big[ 1+\1_{W_n}(y) \big]\,\1_{c_n^{\a}\,r B}(y - x)\,\a^{(2)}({\rm d}(x,y)) \bigg]\phantom{0000000000000000000000000000000000000000000000000000000000000000000000000000000}&&\n\\
\phantom{000000000}&&\n\\
= \frac{\l}{|W_n|^{2\a}} \int\limits_{\R^{3d}}\frac{|W_n \cap (W_n-u)|}{|W_n|}\,\1_{c_n^{\a} r B}(y)
\,\1_{c_n^{\a} r B}(v - u)\,\big(\,\a^{(4)}_{red}({\rm d}(y,u,v)) - \a^{(2)}_{red}({\rm d}y)\,\a^{(2)}({\rm d}(u,v))\,\big)
\phantom{00000000000000000000000000000000000000000000}&&\n\\
\phantom{000000000}&&\n\\
+ \frac{\l}{|W_n|^{2\a}} \int\limits_{\R^{2d}}\bigg(\,1+\frac{|W_n \cap (W_n-y)|}{|W_n|} + 
\frac{|W_n \cap (W_n-u)|}{|W_n|}+\frac{|(W_n-y) \cap (W_n-u)|}{|W_n|}\,\bigg)\, 
\phantom{000000000000000000000000000000000000000000000000000000}&&\n\\
\phantom{000000000}&&\n\\
\times \1_{c_n^{\a} r B}(y)\,\1_{c_n^{\a} r B}(u)\,\a^{(3)}_{red}({\mathrm d}(y,u))+ \frac{\l}{|W_n|^{2\a}}   
\int\limits_{\R^d}\bigg(\,1+\frac{|W_n \cap (W_n - y)|}{|W_n|}\,\bigg)\,\1_{c_n^{\a} r B}(x)\,
\a^{(2)}_{red}({\rm d}x)\phantom{00000000000000000000000000000000000000000000000000}&&
\n\eea

\bea  
\longn 
\left\{\begin{array}{ll}\tau_B^2(r) & \quad\mbox{for}\quad \a = 0\;,\\ 
\phantom{000} &  \label{Kwidetilde}\\
4\,\l^3\,\big(\,1 + \g^{(2)}_{red}(\R^d)\,\big)\,|\,r\,B\,|^2 = 4\,\l^2\,\sigma^2\,|\,r\,B\,|^2 & 
\quad\mbox{for}\quad \a > 0\;.
\end{array}\right.
\eea

\noindent
Both limits in (\ref{Kwidetilde}) follow from (\ref{3red}) and (\ref{a4a2a2})  
combined with a multiple use of  Lemma 1 (for $\a=0$ and $\a > 0$ separately). 
Comparing (\ref{Kwidetilde}) with (\ref{varKdach}) and (\ref{asvar}) reveals that  

\mn
\be    
\lim_{n \to \infty}\frac{|W_n|}{|W_n|^{2\a}}\;\var\big[(\widetilde{\l^2 K_B})_n(c_n^{\a}\,r)\big] 
= \lim_{n \to \infty}\frac{|W_n|}{|W_n|^{2\a}}\;\var\big[(\widehat{\l^2 K_B})_n(c_n^{\a}\,r)\big]\;\;\mbox{for}\;\;0 \le \a < 1\,.
\label{tildehat}\ee

\bigskip
From (\ref{cgh}) and (\ref{a4a2a2}) for $g(x,y)= h(x,y)= |W_n|^{-1}\,\1_{W_n}(x)\,\1_{W_n}(y)\,\1_{c_n^{\a}\,r\,B}(y - x)$  we get in the same way as before
\bea
&&\;\phantom{00000000000000000}|W_n|^{1-2\a}\;\var\bigl[(\overline{\l^2 K_B})_n(c_n^{\a}\,r)\bigr] \n\\
&&\phantom{000000}\n\\
&=&\frac{1}{|W_n|^{1+2\a}}\,\bigg[\,\int\limits_{W_n^4}\1_{c_n^{\a}\,r\,B}(y-x)\,\1_{c_n^{\a}\,r\,B}(v-u)\,
\big(\,\a^{(4)} - \a^{(2)} \times \a^{(2)}\,\big)({\rm d}(x,y,u,v))\n\\
&&\phantom{0000}\n\\
&+& 4\,\int\limits_{W_n^3} \1_{c_n^{\a}\,r B}(y - x)\,\1_{c_n^{\a} r B}(u - x)\,\a^{(3)}({\mathrm d}(x,y,u)) 
+ 2\,\int\limits_{W_n^2} \1_{c_n^{\a} r B}(y - x)\,\a^{(2)}({\rm d}(x,y))\,\bigg] 
\n\eea

\bea
&=&\frac{\l}{|W_n|^{2\a}} \int\limits_{\R^{3d}}\frac{|W_n \cap (W_n-y) \cap (W_n-u) \cap (W_n-v)|}{|W_n|}\,\1_{c_n^{\a} r B}(y) \,\1_{c_n^{\a} r B}(v - u)\,\n\\
&\times& \big(\,\a^{(4)}_{red}({\rm d}(y,u,v)) - \a^{(2)}_{red}({\rm d}y)\,\a^{(2)}({\rm d}(u,v))\,\big)\n\\
&&\phantom{0000}\n\\
&+&\frac{4\,\l}{|W_n|^{2\a}} \int\limits_{\R^{2d}} \frac{|W_n \cap (W_n-y) \cap (W_n-u)|}{|W_n|} \,
\1_{c_n^{\a} r B}(y)\, \1_{c_n^{\a} r B}(u)\,\a^{(3)}_{red}({\mathrm d}(y,u))  \n\\
&&\phantom{0000}\n\\
&+& \frac{2\,\l}{|W_n|^{2\a}}\int\limits_{\R^d}\frac{|W_n \cap (W_n - x)|}{|W_n|}\,
\1_{c_n^{\a} r B}(x)\,\a^{(2)}_{red}({\rm d}x) \;\longn 
\left\{\begin{array}{ll}\tau_B^2(r) & \quad\mbox{for}\quad \a = 0\;,\\ 
\phantom{000} &  \n\\
4\,\l^2\,\sigma^2\,|\,r\,B\,|^2 & \quad\mbox{for}\quad \a > 0\;.
\end{array}\right.
\n\eea

\mn
Hence,
\be    
\lim_{n \to \infty} \frac{|W_n|}{|W_n|^{2\a}}\;\var\big[(\overline{\l^2 K_B})_n(c_n^{\a}\,r)\big] 
= \lim_{n \to \infty}\frac{|W_n|}{|W_n|^{2\a}}\;\var\big[(\widetilde{\l^2 K_B})_n(c_n^{\a}\,r)\big]\;\;\mbox{for}\;\;0 \le \a < 1\,.
\label{linetilde}\ee

\bn
Next we regard the asymptotic behaviour of the covariance of $(\widehat{\l^2 K_B})_n(c_n^{\a}\,r)$ and $(\widetilde{\l^2 K_B})_n(c_n^{\a}\,r)$. For doing this we use (\ref{cgh}) with $g(x,y) = \1_{W_n}(x)\,\1_{W_n}(y)\,\1_{c_n^{\a} r B}(y - x) |W_n \cap (W_n - y + x)|^{-1}$ and  $h(u,v)= |W_n|^{-1} \1_{W_n}(u)\,\1_{c_n^{\a} r B}(v - u)$ which yields

\bea
&&\phantom{00000000000000000}|W_n|^{1-2\a}\;\cov\bigl((\widehat{\l^2 K_B})_n(c_n^{\a}\,r),(\widetilde{\l^2 K_B})_n(c_n^{\a}\,r)\bigr) \n\\
&&\phantom{000000}\n\\
&=&\frac{1}{|W_n|^{2\a}}\int\limits_{\R^{4d}}\frac{\1_{W_n}(x) \1_{W_n}(y) \1_{c_n^{\a} r B}(y-x)}{|W_n \cap (W_n - y + x)|} \1_{W_n}(u)
\1_{c_n^{\a} r B}(v-u) \big(\a^{(4)}-\a^{(2)}\times\a^{(2)}\big)({\rm d}(x,y,u,v))\phantom{00000}\n\\
&&\phantom{0000}\n\\
&+&\frac{2}{|W_n|^{2\a}} \int\limits_{\R^{3d}} \frac{\1_{W_n}(x)\,\1_{W_n}(y)\,\1_{c_n^{\a} r B}(y - x)}{|W_n \cap (W_n - y + x)|}\,\big(\,1+ \1_{W_n}(u)\,\big)\,\1_{c_n^{\a} r B}(u-x)\,\a^{(3)}({\mathrm d}(x,y,u))\n\\
&&\phantom{0000}\n\\
&+&\frac{2}{|W_n|^{2\a}}\int\limits_{\R^{2d}} \frac{\1_{W_n}(x)\,\1_{W_n}(y)\,\1_{c_n^{\a} r B}(y - x)}{|W_n \cap (W_n - y + x)|}\,\a^{(2)}({\rm d}(x,y))
\n\eea

\bea
&=&\frac{\l}{|W_n|^{2\a}}\int\limits_{\R^{3d}}\frac{|W_n \cap (W_n-y) \cap (W_n-u)|}{|W_n \cap (W_n-y)|}\,\1_{c_n^{\a} r B}(y) \,\1_{c_n^{\a} r B}(v - u)\n\\
&\times& \big[\,\a^{(4)}_{red}({\rm d}(y,u,v)) - \a^{(2)}_{red}({\rm d}y)\,\a^{(2)}\big)({\rm d}(u,v))\,\big]\n\\
&&\phantom{0000}\n\\
&+&\frac{2\,\l}{|W_n|^{2\a}} \int\limits_{\R^{2d}}\bigg( 1 + \frac{|W_n \cap (W_n-y) \cap (W_n-u)|}{|W_n \cap (W_n-y)|} \bigg) \,
\1_{c_n^{\a} r B}(y)\,\1_{c_n^{\a} r B}(u)\,\a^{(3)}_{red}({\mathrm d}(y,u))\n\eea
\bea
+ \quad\frac{2\,\l\,\a^{(2)}_{red}(c_n^{\a} r B)}{|W_n|^{2\a}}\;\longn 
\left\{\begin{array}{ll}\tau_B^2(r) & \;\mbox{for}\quad \a = 0\;,\\ 
& \label{covhattilde}\\
4\,\l^2\,\sigma^2\,|\,r\,B\,|^2 & \;\mbox{for}\quad \a > 0\;.
\end{array}\right.\phantom{0000000000000000000000000}
\eea

\bn
To obtain the limit (\ref{covhattilde}) we have used quite the same arguments as in the proof of (\ref{Kwidetilde}). 
Once again we employ (\ref{cgh}) with
$g(x,y)= |W_n|^{-1}\,\1_{W_n}(x)\,\1_{W_n}(y)\,\1_{c_n^{\a} r B}(y - x) $ and  $h(u,v)= |W_n|^{-1}\,\1_{W_n}(u)\,\1_{c_n^{\a} r B}(v - u)$. 
After performing almost the same calculations as before we arrive at 

\bea
&&\phantom{00000000000000000}|W_n|^{1-2\a}\;\cov\bigl((\overline{\l^2 K_B})_n(c_n^{\a}\,r),(\widetilde{\l^2 K_B})_n(c_n^{\a}\,r)\bigr) \n\\
&&\phantom{000000}\n\\
&=&\frac{1}{|W_n|^{2\a}}\int\limits_{\R^{4d}}\frac{\1_{W_n}(x)\1_{W_n}(y)\1_{c_n^{\a}r B}(y-x)}{|W_n|}\1_{W_n}(u)
\1_{c_n^{\a} r B}(v-u)\big(\a^{(4)} - \a^{(2)}\times\a^{(2)}\big)({\rm d}(x,y,u,v))\phantom{0000}\n\\
&&\phantom{0000}\n\\
&+&\frac{2}{|W_n|^{2\a}} \int\limits_{\R^{3d}} \frac{\1_{W_n}(x)\,\1_{W_n}(y)\,\1_{c_n^{\a} r B}(y - x)}{|W_n|}\,\big(\,1+ \1_{W_n}(u)\,\big)\,\1_{c_n^{\a} r B}(u-x)\,\a^{(3)}({\mathrm d}(x,y,u))\n\\
&&\phantom{0000}\n\\
&+&\frac{2}{|W_n|^{2\a}}\int\limits_{\R^{2d}} \frac{\1_{W_n}(x)\,\1_{W_n}(y)\,\1_{c_n^{\a} r B}(y - x)}{|W_n|}\,\a^{(2)}({\rm d}(x,y))\n\\
&&\phantom{0000}\n\\
&=&\frac{\l}{|W_n|^{2\a}}\int\limits_{\R^{3d}}\frac{|W_n \cap (W_n-y) \cap (W_n-u)|}{|W_n|}\,\1_{c_n^{\a} r B}(y) \,\1_{c_n^{\a} r B}(v - u)\n\\
&\times& \big[\,\a^{(4)}_{red}({\rm d}(y,u,v)) - \a^{(2)}_{red}({\rm d}y)\,\a^{(2)}\big)({\rm d}(u,v))\,\big]\n\\
&&\phantom{0000}\n\\
&+&\frac{2\,\l}{|W_n|^{2\a}} \int\limits_{\R^{2d}}\bigg( 1 + \frac{|W_n \cap (W_n-y) \cap (W_n-u)|}{|W_n|} \bigg) \,
\1_{c_n^{\a} r B}(y)\,\1_{c_n^{\a} r B}(u)\,\a^{(3)}_{red}({\mathrm d}(y,u))\n\eea
\bea
+\quad\frac{2\,\l}{|W_n|^{2\a}}\int\limits_{c_n^{\a} r B}\frac{| W_n \cap (W_n -y)|}{|W_n|}\,\a^{(2)}({\rm d}(x,y))\;\longn 
\left\{\begin{array}{ll}\tau_B^2(r) & \;\mbox{for}\quad \a = 0\;,\\ 
\phantom{000} & \label{covlinetilde}\\
4\,\l^2\,\sigma^2\,|\,r\,B\,|^2 & \;\mbox{for}\quad \a > 0\;.
\end{array}\right.\phantom{00}
\eea

\mn
Since the limits (\ref{covhattilde}) and (\ref{covlinetilde}) are identical and coincide with the asymptotic  variances 
(\ref{asvar}) and (\ref{tildehat}) for $0 \le \a < 1$  the first assertion of Lemma 2B is proved due to the well-known relation 
\bea 
\ew\bigl[(\widehat{\l^2 K_B})_n(c_n^{\a}\,r) - (\widetilde{\l^2 K_B})_n(c_n^{\a}\,r)\bigr]^2 &=& \var\big[(\widehat{\l^2 K_B})_n(c_n^{\a}\,r)\big]  
+ \var\big[(\widetilde{\l^2 K_B})_n(c_n^{\a}\,r)\big]\n\\
&-& 2\,\cov\bigl((\widehat{\l^2 K_B})_n(c_n^{\a}\,r),(\widetilde{\l^2 K_B})_n(c_n^{\a}\,r)\bigr)\,.
\n\eea
In just the same way the identity
\bea
\var\bigl[(\overline{\l^2 K_B})_n(c_n^{\a}\,r) - (\widetilde{\l^2 K_B})_n(c_n^{\a}\,r)\bigr]^2 &=& \var\big[(\overline{\l^2 K_B})_n(c_n^{\a}\,r)\big]  
+ \var\big[(\widetilde{\l^2 K_B})_n(c_n^{\a}\,r)\big]\n\\
&-& 2\,\cov\bigl((\overline{\l^2 K_B})_n(c_n^{\a}\,r),(\widetilde{\l^2 K_B})_n(c_n^{\a}\,r)\bigr)
\n\eea
implies the second assertion of Lemma 1B which completes the proof of Lemma 1B. $\Box$

\bn
\textbf{Lemma 2B} \ Under the conditions of Lemma 2 we have
\be
\lim_{n \to \infty}\,|W_n|^{1 - 2\a}\,\ew\big(\,(\widehat{\l^2 K_B})_n(c_n^{\a}\,r)\, -
\,(\widehat{\l^2})_n\,K_B(c_n^{\a}\,r)\,\big)^2\,=\,0 \quad\mbox{for all}\quad r \ge 0\,.
\label{hatlambda2}\ee
Furthermore, $(\widehat{\l^2 K_B})_n(c_n^{\a} r)$ can be replaced by $(\widetilde{\l^2 K_B})_n(c_n^{\a} r)\,$.

\bn
{\em Proof of Lemma 2B} \ By  $(a + b + c)^2 \le 3\,a^2 + 3\,b^2 +3\,c^2$  and $(\widehat{\l^2})_n = (\widehat{\l}_n)^2 - \widehat{\l}_n/|W_n|\,$,

\bea 
&&\qquad \;|W_n|^{1 - 2\a}\;\ew\big(\,(\widehat{\l^2 K_B})_n(c_n^{\a}\,r) - (\widehat{\l^2})_n\,K_B(c_n^{\a}\,r)\,\big)^2\n\\
&&\phantom{00}\n\\
&=& |W_n|^{1-2\a}\,\ew\Big(\,(\widehat{\l^2 K_B})_n(c_n^{\a}\,r) - \l^2\,K_B(c^{\a}_n\,r) - 
\big((\widehat{\l^2})_n - \l^2 \big)\,K_B(c_n^{\a}\,r)\,\Big)^2\n\\
&&\phantom{00}\n\\
&=&  |W_n|^{1-2\a}\;\ew\Big(\,(\widehat{\l^2 K_B})_n(c_n^{\a}\,r) - \l^2\,K_B(c^{\a}_n\,r) - 2\,\l\,K_B(c^{\a}_n\,r)\,\big( \widehat{\l}_n - \l \big)\n\\
&&\qquad - \,K_B(c^{\a}_n\,r)\,\big((\widehat{\l^2})_n - \l^2 - 2\,\l\,(\widehat{\l}_n - \l)\big)\,\Big)^2\n\\
&&\phantom{00}\n\\
&\le& 3\,|W_n|^{1-2\a}\;\ew\Big(\,(\widehat{\l^2 K_B})_n(c_n^{\a}\,r) - \l^2\,K_B(c^{\a}_n\,r) - 2\,\l\,K_B(c^{\a}_n\,r)\,\big( \widehat{\l}_n - \l \big)\,\Big)^2\n\\
&&\qquad + \,|W_n|^{1-2\a}\;\ew\big(\,K_B(c^{\a}_n\,r)\,\big(\widehat{\l}_n - \l\big)^2\,\big)^2 + 
3\,|W_n|^{-(1+2\a)}\;\ew\big(\,K_B(c^{\a}_n\,r)\,\widehat{\l}_n)\,\big)^2\,.
\label{lambda2}\eea
 
\bigskip
Clearly, $|W_n|\,\ew \big(\widehat{\l}_n - \l\big)^4 = |W_n|^{-3}\,\ew \big(\,N(W_n) -\l\,|W_n|\,\big)^4$
and it is easily checked that $\ew \big(\,N(W_n) -\l\,|W_n|\,\big)^4 = \mathbf{Cum}_4(N(W_n)) + 3\,\big(\,\var(N(W_n))\,\big)^2$, where 
$\mathbf{Cum}_4(X) = \frac{{\rm d}^4}{{\rm d}t^4} \log \ew e^{it X}|_{t=0} $ denotes the fourth cumulant of the random variable $X$. 
The fourth cumulant of $N(W_n)$  can be expressed as linear combination of the  factorial cumulant measures $\g^{(k)}$ applied to $k$-fold Cartesian
product $W_n^k$  and $\ew N(W_n)$. A somewhat lengthy calculation based on (\ref{ga}) shows  that
\be 
\mathbf{Cum}_4(N(W_n)) = \g^{(4)}(W_n\times\cdots\times W_n) + 6\,\g^{(3)}(W_n\times W_n\times W_n) + 7\,\g^{(2)}(W_n\times W_n) + \ew N(W_n)\,.
\n\ee
The $B_4$-mixing property of the PP   $N = \sum_{i\ge 1}\d_{X_i}$ implies immediately the estimates
\be 
|\,\mathbf{Cum}_4(N(W_n))\,| \le \l\,|W_n|\,\big(\,\g_{red}^{(4)}((\R^d)^3) + 6\,\g_{red}^{(3)}((\R^d)^2) + 7\,\g_{red}^{(2)}(\R^d) + 1\,\big)
\n\ee
and $\var(N(W_n)) = \g_2(W_n\times W_n) + \ew N(W_n) \le \l\,|W_n|\,\big(\,\g_{red}^{(2)}(\R^d) + 1\,\big)$. These two estimates combined with (\ref{ared}) yield
\be
|W_n|^{1-2\a}\,\ew\big(\,K_B(c^{\a}_n\,r)\,\big(\widehat{\l}_n - \l\big)^2\,\big)^2 =
\Big(\frac{K_B(c^{\a}_n\,r)}{|W_n|^{\a}}\Big)^2\,|W_n|\,\ew\big(\widehat{\l}_n - \l\big)^4 = {\mathcal O}(|W_n|^{-1})
\n\ee
and
\be
|W_n|^{-(1+2\a)}\;\ew\big(\,K_B(c^{\a}_n\,r)\,\widehat{\l}_n)\,\big)^2 = \Big(\frac{K_B(c^{\a}_n\,r)}{|W_n|^{\a}}\Big)^2\,\frac{\ew(N(W_n))^2}{|W_n|^3}
 =  {\mathcal O}(|W_n|^{-1})
\n\ee 

\mn
as $n \to \infty$. Lemma 2 and the latter two relations show that the r.h.s. of (\ref{lambda2})
converges to zero which terminates the proof (\ref{hatlambda2}). Finally, by the norm inequality in $L^2$ it follows

\bea 
\Big(\ew\big(\,(\widetilde{\l^2 K_B})_n(c_n^{\a}\,r) - (\widehat{\l^2})_n\,K_B(c_n^{\a}\,r)\,\big)^2\Big)^{1/2}
&\le & \Big(\ew\big(\,(\widehat{\l^2 K_B})_n(c_n^{\a}\,r) - (\widetilde{\l^2 K_B})_n(c_n^{\a}\,r)\,\big)^2\Big)^{1/2}\n\\
&+& \Big(\ew\big(\,(\widehat{\l^2 K_B})_n(c_n^{\a}\,r) - (\widehat{\l^2})_n\,K_B(c_n^{\a}\,r)\,\big)^2 \Big)^{1/2}\,.
\n\eea
After multiplying both sides by $|W_n|^{1/2+\a}$ the first limit of Lemma 1B and (\ref{hatlambda2}) complete the proof of Lemma 2B. $\Box$  

\bigskip
At the end of this Appendix B we given a structural formula for the covariance function of the Gaussian limit process of 
$\sqrt{|W_n|}\,\big((\widehat{\l^2 K_B})_n(r) - \l^2 K_B(r)\big)$ for $r \ge 0\,$, which can be obtained simply by repeating 
the proof of (\ref{asvar}) taking into account some obvious changes.

\mn
\textbf{Lemma 3B} \ For any  $B_4$-mixing PP  $N = \sum_{i\ge 1} \d_{X_i}$ on $\R^d$ the asymptotic covariance 
\be 
\tau_B(s,t) := \lim_{n \to \infty} |W_n|\,\cov\bigl(\,(\widehat{\l^2 K_B})_n(s),(\widehat{\l^2 K_B})_n(t)\,\bigr) 
\n\ee
exists for any $s,t \ge 0$ and takes the form

\bea 
\tau_B(s,t) & = &
\l\,\int_{\R^{3\,d}}\1_{s\,B}(y-x)\,\1_{t\,B}(z)\,\g^{(4)}_{red}({\rm d}(x,y,z)) + 2\,\l^2\,|s B|\,\g^{(3)}_{red}(t B \times \R^d)\n\\
&+& 2\,\l^2\,|t B|\,\g^{(3)}_{red}(s B \times \R^d)
+ 2\,\l^3\,\int_{\R^{2\,d}}|(s\,B-x) \cap (t\,B-y) |\,\g^{(2)}_{red}({\rm d}x)\,\g^{(2)}_{red}({\rm d}y)\n\\
&+& 4\,\l^3\,|s\,B|\,|t\,B|\,\g^{(2)}_{red}(\R^d) + 4\,\l\,\a^{(3)}_{red}(s\,B \times t\,B)
+ 2\,\l\,\a^{(2)}_{red}((s \wedge t)\,B)\,.
\label{ascov}\eea  

\medskip
We notice that $\int_{\R^{3\,d}}\1_{s\,B}(y-x)\,\1_{t\,B}(z)\,\g^{(4)}_{red}({\rm d}(x,y,z)) = \int_{\R^{3\,d}}\1_{t\,B}(y-x)\,\1_{s\,B}(z)\,\g^{(4)}_{red}({\rm d}(x,y,z))$ 
for all $s,t \ge 0$ due to the symmetry properties of $\g^{(4)}_{red}(\cdot)\,$.

\bigskip
In case of stationary Poisson cluster PPes and ($\a$-)determinantal PPes the reduced cumulant measures $\g_{red}^{(2)}, \g_{red}^{(3)}, \g_{red}^{(4)}$ 
resp. their Lebesgue densities have a comparatively simple shape leading to more compact  representations of $\tau_B(s,t)$, see \cite{he88}, \cite{he16} or 
\cite{bl16, bl17}. In the special case of a stationary Poisson PP with intensity $\l > 0$  (implying $\g_{red}^{(k)}(\cdot) \equiv 0$ for $k=2,3,4$ ) we get 
\be
\tau_B(s,t) = 2\,\l^2\,(s \wedge t)^d\,| B |\,\big(\,1 + 2\,\l\,(s \vee t)^d\,| B |\,\big)\quad\mbox{for}\quad s,t \ge 0\;\;,\;\; \mbox{see \,\cite{he91}\,,\,\cite{he13}}\,.
\n\ee

\vfill\eject

\end{document}